%% file: Stevens-semistable.tex
\documentclass[11pt, reqno]{amsart}

\usepackage[OT2, T1]{fontenc}

\usepackage[usenames,dvipsnames]{color}

\usepackage{url}
\usepackage{amsmath}
\usepackage{array}
\usepackage{graphicx}
\usepackage{amssymb}
\usepackage{amsthm}
\definecolor{link}{RGB}{11,0,128}
\usepackage[colorlinks=true,citecolor=NavyBlue,linkcolor=Bittersweet,urlcolor=link]{hyperref}
\usepackage{paralist}
\usepackage{colonequals}		
\usepackage{color}
\usepackage{mathabx}		
\usepackage{amscd}
\usepackage[all,cmtip]{xy}	
\usepackage{sseq}			
\usepackage{verbatim}
\usepackage{parskip}
\usepackage{microtype}
\usepackage[in]{fullpage}
\usepackage{mathrsfs} 			
\usepackage{cleveref}
\usepackage{enumitem}
\usepackage{moreenum}			
\usepackage[alphabetic,lite]{amsrefs} 	
\usepackage{stmaryrd}	
\usepackage[foot]{amsaddr}
\usepackage{subfiles}
\usepackage{ragged2e}

\DeclareSymbolFont{cyrletters}{OT2}{wncyr}{m}{n}
\DeclareMathSymbol{\Sha}{\mathalpha}{cyrletters}{"58}

\newcommand{\bC}{\mathbb{C}}

\newcommand{\bF}{\mathbb{F}}
\newcommand{\bG}{\mathbb{G}}

\newcommand{\bP}{\mathbb{P}}
\newcommand{\bQ}{\mathbb{Q}}

\newcommand{\bZ}{\mathbb{Z}}

\newcommand{\cE}{\mathcal{E}}

\newcommand{\cJ}{\mathcal{J}}

\newcommand{\sH}{\mathscr{H}}

\newcommand{\sO}{\mathscr{O}}

\newcommand{\sS}{\mathscr{S}}

\newcommand{\sU}{\mathscr{U}}
\newcommand{\sV}{\mathscr{V}}

\newcommand{\sX}{\mathscr{X}}
\newcommand{\sY}{\mathscr{Y}}
\newcommand{\sZ}{\mathscr{Z}}

\newcommand{\ra}{\rightarrow}

\newcommand{\xra}{\xrightarrow}

\newcommand{\hra}{\hookrightarrow}

\newcommand{\wt}{\widetilde}
\newcommand{\wh}{\widehat}

\newcommand{\ce}{\colonequals}
\newcommand{\ov}{\overline}

\newcommand{\sm}{\mathrm{sm}}

\renewcommand{\b}{\textbf}
\newcommand{\surjects}{\twoheadrightarrow}

\newcommand{\tensor}{\otimes} 		
\newcommand{\isomto}{\overset{\sim}{\longrightarrow}}

\newcommand{\et}{\mathrm{\acute{e}t}}	
\newcommand{\llb}{\llbracket}		
\newcommand{\rrb}{\rrbracket}		
\newcommand{\sh}{\mathrm{sh}}		

\renewcommand{\i}{^{-1}}
\renewcommand{\th}{^{\mathrm{th}}}

\providecommand{\abs}[1]{\left\lvert#1\right\rvert}

\providecommand{\p}[1]{\left(#1\right)}

\providecommand{\f}[2]{\frac{#1}{#2}}
\DeclareMathOperator{\Ker}{Ker}			
\DeclareMathOperator{\Coker}{Coker}		
\DeclareMathOperator{\im}{Im}			
\DeclareMathOperator{\Spec}{Spec}		
\DeclareMathOperator{\Hom}{Hom}			
\DeclareMathOperator{\ord}{ord}	
\DeclareMathOperator{\GL}{GL}		
\DeclareMathOperator{\Aut}{Aut}		
\DeclareMathOperator{\Lie}{Lie}		
\DeclareMathOperator{\Pic}{Pic}		
\newcommand{\ba}{\begin{aligned}}
\newcommand{\ea}{\end{aligned}}
\newcommand{\be}{\begin{equation}}
\newcommand{\ee}{\end{equation}}
\newcommand{\pf}{\begin{proof}}
\newcommand{\bpf}{\begin{proof}}
\newcommand{\epf}{\end{proof}}
\newcommand{\bthm}{\begin{thm}}
\newcommand{\ethm}{\end{thm}}
\newcommand{\bthmt}{\begin{thm-tweak}}
\newcommand{\ethmt}{\end{thm-tweak}}
\newcommand{\bconjt}{\begin{conj-tweak}}
\newcommand{\econjt}{\end{conj-tweak}}
\newcommand{\bprop}{\begin{prop}}
\newcommand{\eprop}{\end{prop}}
\newcommand{\bcor}{\begin{cor}}
\newcommand{\ecor}{\end{cor}}

\newcommand{\brem}{\begin{rem}}
\newcommand{\erem}{\end{rem}}
\newcommand{\bremt}{\begin{rem-tweak}}
\newcommand{\eremt}{\end{rem-tweak}}
\newcommand{\brems}{\begin{rems} \hfill \begin{enumerate}[label=\b{\thesubsubsection.},ref=\thesubsubsection]}
\newcommand{\bremstweak}{\begin{rems-tweak} \hfill \begin{enumerate}[label=\b{\thesubsection.},ref=\thesubsection]}
\newcommand{\bremst}{\begin{rems-tweak} \hfill \begin{enumerate}[label=\b{\thesubsection.},ref=\thesubsection]}

\newcommand{\remit}{\addtocounter{subsection}{1} \item}
\newcommand{\erems}{\end{enumerate} \end{rems}}
\newcommand{\eremstweak}{\end{enumerate} \end{rems-tweak}}
\newcommand{\eremst}{\end{enumerate} \end{rems-tweak}}
\newcommand{\blem}{\begin{lemma}}
\newcommand{\elem}{\end{lemma}}
\newcommand{\blemt}{\begin{lemma-tweak}}
\newcommand{\elemt}{\end{lemma-tweak}}
\newcommand{\bconj}{\begin{conj}}
\newcommand{\econj}{\end{conj}}
\newcommand{\bprob}{\begin{Problem}}
\newcommand{\eprob}{\end{Problem}}
\newcommand{\bpropt}{\begin{prop-tweak}}
\newcommand{\epropt}{\end{prop-tweak}}
\newcommand{\bq}{\begin{q}}
\newcommand{\eq}{\end{q}}
\newcommand{\benum}{\begin{enumerate}[label={{\upshape(\alph*)}}]}
\newcommand{\benuma}{\begin{enumerate}[label={{\upshape(\arabic*)}}]}
\newcommand{\benumr}{\begin{enumerate}[label={{\upshape(\roman*)}}]}
\newcommand{\eenum}{\end{enumerate}}
\providecommand{\up}[1]{{\upshape(}#1{\upshape)}}
\providecommand{\uref}[1]{{\upshape\ref{#1}}}
\providecommand{\uS}{{\upshape\S}}
\providecommand{\ucolon}{{\upshape:} }
\providecommand{\uscolon}{{\upshape;} }
\newcommand{\bc}{}
\newcommand{\bd}{\begin{defn}}
\newcommand{\ed}{\end{defn}}

\newcommand{\beg}{\begin{eg}}
\newcommand{\eeg}{\end{eg}}
\newcommand{\begt}{\begin{eg-tweak}}
\newcommand{\eegt}{\end{eg-tweak}}
\newcommand{\bcl}{\begin{claim}}
\newcommand{\ecl}{\end{claim}}
\newcommand{\lab}{\label}

\newcommand{\q}{\quad}
\newcommand{\qq}{\quad\quad}
\newcommand{\qqq}{\quad\quad\quad}

\newcommand{\tst}{\textstyle}

\newcommand{\forg}{\mathrm{forg}}
\newcommand{\quot}{\mathrm{quot}}

\newcommand{\cusps}{\mathrm{cusps}}

\newcommand*{\QED}{\hfill\ensuremath{\qed}}
\newcommand*{\QEDD}{\hfill\ensuremath{\qed\qed}}

\theoremstyle{plain}
\newtheorem{thm}[subsubsection]{Theorem}
\Crefname{thm}{Theorem}{Theorems}

\newtheorem{thm-tweak}[subsection]{Theorem}
\Crefname{thm-tweak}{Theorem}{Theorems}

\Crefname{rethm}{Theorem}{Theorem}
\newtheorem{prop}[subsubsection]{Proposition}
\Crefname{prop}{Proposition}{Propositions} 

\newtheorem{prop-tweak}[subsection]{Proposition}
\Crefname{prop-tweak}{Proposition}{Propositions}

\Crefname{Q}{Question}{Questions}
\Crefname{eg}{Example}{Examples}
\newtheorem{Problem}[subsection]{Problem}
\Crefname{Problem}{Problem}{Problems}
\newtheorem{conj}[subsubsection]{Conjecture}
\Crefname{conj}{Conjecture}{Conjectures}
\newtheorem{conj-tweak}[subsection]{Conjecture}
\Crefname{conj-tweak}{Conjecture}{Conjectures}
\newtheorem{cor}[subsubsection]{Corollary}
\Crefname{cor}{Corollary}{Corollaries}
\newtheorem{cor-tweak}[subsection]{Corollary}
\Crefname{cor-tweak}{Corollary}{Corollaries}

\newtheorem{lemma}[subsubsection]{Lemma}

\newtheorem{lemma-tweak}[subsection]{Lemma}
\Crefname{lemma-tweak}{Lemma}{Lemmas}

\Crefname{subprop}{Proposition}{Propositions}
\Crefname{subcor}{Corollary}{Corollaries}
\Crefname{sublem}{Lemma}{Lemmas}

\theoremstyle{remark}
\newtheorem{claim}[equation]{Claim}
\Crefname{claim}{Claim}{Claims}

\Crefname{subrem}{Remark}{Remarks}

\theoremstyle{definition}
\newtheorem{defn}[subsubsection]{Definition}
\Crefname{defn}{Definition}{Definitions}
\newtheorem{defn-tweak}[subsection]{Definition}
\Crefname{defn-tweak}{Definition}{Definitions}

\Crefname{conv}{Convention}{Conventions}
\newtheorem{eg}[subsubsection]{Example}
\newtheorem{rem}[subsubsection]{Remark}
\Crefname{rem}{Remark}{Remarks}

\newtheorem{eg-tweak}[subsection]{Example}
\Crefname{eg-tweak}{Example}{Examples}

\newtheorem{rem-tweak}[subsection]{Remark}
\Crefname{rem-tweak}{Remark}{Remarks}

\newtheorem*{rems}{Remarks}

\newtheorem*{rems-tweak}{Remarks}

\newtheoremstyle{subsection-tweak}
   {11pt}
   {3pt}%
   {}
   {}%
   {\bfseries}
   {}%
   {.5em}
   {\thmnumber{\@{#1}{}\@{#2}.}%
    \thmnote{~{\bfseries#3.}}}

\Crefname{innercustomconj}{Conjecture}{Conjecture}

\theoremstyle{subsection-tweak}
\newtheorem{pp}[subsection]{}
\newcommand{\bpp}{\begin{pp}}
\newcommand{\epp}{\end{pp}}

\theoremstyle{subsection-tweak}
\newtheorem{pp-tweak}[subsection]{}

\numberwithin{equation}{subsection}

\makeatletter
\def\@tocline#1#2#3#4#5#6#7{
    \begingroup 
    \@ifempty{#4}{%
    }{%
    }%

    \parindent\z@ \leftskip#3\relax \advance\leftskip\@tempdima\relax
    #5\hskip-\@tempdima
      \ifcase #1
       \or\or \hskip 2em \or \hskip 1em \else \hskip 3em \fi%
      #6\nobreak\relax
    \dotfill\hbox to\@pnumwidth{\@tocpagenum{#7}}\par
    \nobreak
    \endgroup
  }
 \def\l@section{\@tocline{1}{0pt}{1pc}{}{}}

\renewcommand{\tocsection}[3]{%
  \indentlabel{\@ifnotempty{#2}{\makebox[1.3em][l]{%
    \ignorespaces#1 \bfseries{#2}.\hfill}}}\bfseries{#3}
    \vspace{1.5pt}}

\renewcommand{\tocsubsection}[3]{%
  \indentlabel{\@ifnotempty{#2}{\hspace*{-0.5em}\makebox[2.1em][l]{%
    \ignorespaces#1#2.\hfill}}}#3
    \vspace{1.5pt}}

\makeatother

\begin{document}
\author{K\k{e}stutis \v{C}esnavi\v{c}ius}
\title{The Manin--Stevens constant in the semistable case}
\date{\today}
\address{Centre national de la recherche scientifique, Institut de Math\'ematique d'Orsay, F-91400,  France}
\email{kestutis@math.u-psud.fr}

\begin{abstract} 
Stevens conjectured that for every optimal parametrization $\phi\colon X_1(n) \surjects E$ of an elliptic curve $E$ over $\bQ$ of conductor $n$, the pullback of some N\'{e}ron differential on $E$ is the differential associated to the normalized new eigenform that corresponds to the isogeny class of $E$. We prove this conjecture under the assumption that $E$ is semistable, the key novelty lying in the $2$-primary analysis when $n$ is even. For this analysis, we first relate the general case of the conjecture to a divisibility relation between $\deg \phi$ and a certain congruence number and then reduce the semistable case to a question of exhibiting enough suitably constrained oldforms. Our methods also apply to parametrizations by $X_0(n)$ and prove new cases of the Manin conjecture.
\end{abstract}

\maketitle

\section{Introduction}

With the purpose of relating the arithmetic of an elliptic curve $E$ over $\bQ$ of conductor $n$ to the arithmetic of the modular curve $X_1(n)$ via a given modular parametrization
\[
\phi\colon X_1(n) \surjects E,
\]
one normalizes by arranging that the cusp $\text{``}0\text{''} \in X_1(n)(\bQ)$ maps to $0 \in E(\bQ)$ and, at the cost of replacing $E$ by an isogenous curve, that the induced quotient map 
\[
\pi\colon J_1(n) \surjects E
\]
from the Jacobian is \emph{optimal} in the sense that its kernel is connected. For such a $\phi$ (equivalently,~$\pi$), one seeks to understand the differential aspect of the modularity relationship captured by the equality
\be \lab{omega-f} \tag{$\bigstar$}
\pi^*(\omega_\cE) = c_\pi \cdot f_E \qq \text{for some}\qq c_\pi \in \bQ^\times,
\ee
where $\omega_\cE \in H^0(E, \Omega^1)$ is a N\'{e}ron differential and $f_E \in H^0(X_1(n), \Omega^1) \cong H^0(J_1(n), \Omega^1)$ is the differential form associated to the normalized new eigenform that corresponds to the isogeny class of $E$. Since $\pi$ is \emph{new}, that is, factors through the new quotient of $J_1(n)$, the multiplicity one principle supplies \eqref{omega-f} and it remains to understand the appearing \emph{Manin--Stevens constant} $c_\pi$.

\bconjt[Stevens, \cite{Ste89}*{Conj.~I$'$~(a)}] \lab{stevens-conj}
For a new elliptic optimal quotient $\pi\colon J_1(n) \surjects E$,
\[
c_\pi = \pm 1.
\]
\econjt


We settle \Cref{stevens-conj} for semistable $E$ and, more generally, settle its $p$-primary part for primes $p$ with $\ord_p(n) \le 1$. For this, existing techniques suffice if $p$ is odd, so the key new case is $p = 2$.

\bthmt[\S\ref{stevens-pf-1} and \S\ref{stevens-pf-2}] \lab{stevens-main}
For a new elliptic optimal quotient $\pi\colon J_1(n) \surjects E$ and a prime $p$, 
\[
\text{if} \q \ord_p(n) \le 1, \q \text{then} \q \ord_p(c_\pi) = 0.
\]
In particular, if $E$ is semistable \up{that is, if $n$ is squarefree}, then $c_\pi = \pm 1$.
\ethmt

\Cref{stevens-conj} of Stevens is a variant of an earlier conjecture of Manin on parametrizations by $X_0(n)$ (equivalently, by its Jacobian $J_0(n)$), for which the analogous $c_\pi \in \bQ^\times$  is the \emph{Manin constant}.

\bconjt[Manin, \cite{Man71}*{10.3}] \lab{manin-conj}
For a new elliptic optimal quotient $\pi\colon J_0(n) \surjects E$, 
\[
c_\pi = \pm 1.
\]
\econjt

One knows that \Cref{manin-conj} implies \Cref{stevens-conj} (see \Cref{jj-ee}~\ref{JE-b}). Cremona has proved \Cref{manin-conj} for all $E$ of conductor at most $130000$, see \cite{ARS06}*{Thm.~2.6} (see also \cite{Cre16} for the same result up to conductor $370000$).

One typically studies the Manin constant through its $p$-adic valuations. The following theorem summarizes the known cases of the $p$-primary part of \Cref{manin-conj} at the semistable primes $p$.

\bthmt \lab{manin-known}
For a new elliptic optimal quotient $\pi\colon J_0(n) \surjects E$ and a prime $p$,
\[
\text{if} \q \ord_p(n) \le 1, \q  \text{then} \q \ord_p(c_\pi) = 0 \q \text{in any of the following cases\ucolon}
\]
\benumr
\item \lab{MK-1} {\upshape (Mazur, \cite{Maz78}*{Cor.~4.1}; see also \Cref{odd-reproof} and \Cref{odd-p}).}
If $p$ is odd\uscolon

\item \lab{MK-2} {\upshape (Abbes--Ullmo, \cite{AU96}*{Thm.~A}; see also \Cref{odd-n}).}
If $p = 2$ and $\ord_2(n) = 0$\uscolon

\item \lab{MK-3} {\upshape (Raynaud, \cite{AU96}*{(ii) on p.~270}).}
If $p = 2$ and $\ord_2(\Delta_E)$ is odd, where $\Delta_E \in \bQ^\times$ is the discriminant of a Weierstrass equation for $E$\uscolon

\item \lab{MK-4} {\upshape (Agashe--Ribet--Stein, \cite{ARS06}*{Thm.~2.7}; see also Remark \ref{degree-divides}).}
If $p = 2$ and the degree of the composite $X_0(n) \hra J_0(n) \xra{\pi} E$ obtained by choosing a point in $X_0(n)(\bQ)$ is odd.
\eenum
\ethmt

In addition, for a $\pi$ as in \Cref{manin-known}, one knows that $0 \le \ord_2(c_\pi) \le 1$ when $\ord_2(n) \le 1$ thanks to a result of Mazur--Raynaud, \cite{AU96}*{Prop.~3.1}, based on exactness properties of semiabelian N\'{e}ron models. The techniques of the proof of \Cref{stevens-main} reprove this result in \Cref{Raynaud} without using such exactness properties. \Cref{odd-reproof} achieves the same for \Cref{manin-known}~\ref{MK-1}. 

Beyond the semistable $p$, for a $\pi$ as in \Cref{manin-known}, Edixhoven proved in \cite{Edi91}*{Thm.~3} that $\ord_p(c_\pi) = 0$ in the case when $p > 7$ and $E_{\bQ_p}$ does not have potentially ordinary reduction of Kodaira type II, III, or IV. Further cases may be supplied by  the unfinished manuscript \cite{Edi01}.

In addition to streamlined reproofs of \Cref{manin-known}~\ref{MK-1}, \ref{MK-2}, and \ref{MK-4}, our methods also lead to the following new cases of the $2$-primary part of the Manin conjecture.\footnote{The semistable case of the Manin conjecture has now been settled in full in \cite{Manin-semistable} by a different method.}

\bthmt \lab{manin-main}
For a new elliptic optimal quotient $\pi\colon J_0(n) \surjects E$,
\[
\text{if} \q \ord_2(n) \le 1, \q  \text{then} \q \ord_2(c_\pi) = 0 \q \text{in any of the following cases\ucolon}
\]
\benumr
\item \lab{MM-1} {\upshape (\S\ref{manin-pf-1}).}
If $n$ has a prime factor $q$ with $q\equiv 3 \bmod 4$\uscolon

\item \lab{MM-15} {\upshape (\S\ref{manin-pf-1}).}
If $n = 2p$ for some prime $p$\uscolon

\item \lab{MM-2} {\upshape (\S\ref{manin-pf-2}).}
If $E(\bQ)[2] = 0$.
\eenum
\ethmt

For further input which would prove that $\ord_2(c_\pi) = 0$ whenever $\ord_2(n) \le 1$, see Remark \ref{always}.

The conditions \ref{MM-1}--\ref{MM-2} in \Cref{manin-main} are global, so the combination of \Cref{manin-known,manin-main} covers significantly more new elliptic optimal quotients than \Cref{manin-known} alone.

\begt
To get a sense of the scope of \Cref{manin-main} we used the website \cite{LMFDB} to inspect all elliptic curves over $\bQ$ of conductor $\le 200$ that are optimal with respect to $X_0(n)$ and semistable at $2$; we found $205$ such curves. For $143$ of them \Cref{manin-known} proves that $\ord_2(c_\pi) = 0$. \Cref{manin-main}~\ref{MM-1} then proves that $\ord_2(c_\pi) = 0$ for $47$ of the remaining $62$ curves $30.a8,\, 34.a4,\, \ldots,\, 198.d4,\, 198.e3$, leaving $15$ curves: $34.a4,\, 58.a1,\, \ldots,\, 178.b2,\, 194.a2$. \Cref{manin-main}~\ref{MM-15} proves that $\ord_2(c_\pi) = 0$ for $10$ of these, leaving $5$ curves: $130.a2,\, 130.b4,\, 130.c1,\, 170.a2,\, 170.b1$, all of which have $E(\bQ)[2] \neq 0$. 

\Cref{manin-main}~\ref{MM-2} does provide new information for some curves, for instance, for $530.a1$, which has $E(\bQ)[2] = 0$ but for which neither \Cref{manin-known}~\ref{MK-2}--\ref{MK-4} nor \Cref{manin-main}~\ref{MM-1}--\ref{MM-15} apply.
\eegt

\bpp[The overview of the proofs]
The first step of the proofs of \Cref{stevens-main,manin-main} is a reduction, not specific to semistable $p$, to a divisibility relation between $\deg \phi$ and a certain congruence number (see \Cref{congJ-deg}~\ref{CD-c}). However, the required divisibility differs from the ones available in the literature because we measure congruences between weight $2$ cusp forms with respect to the cotangent space at the identity of the N\'{e}ron model of $J_1(n)$ (resp.,~of $J_0(n)$) rather than with respect to $q$-expansions at ``$\infty$.'' The proofs proceed to isolate a module that controls the difference between the two types of congruences and, under a semistability assumption, the problem becomes that of exhibiting its vanishing (see \Cref{main-reduction}). For this, it suffices to show that oldforms offset the difference between two integral structures on the $\bQ$-vector space of weight $2$ cusp forms (see the introduction of \S\ref{oldforms}). The technical heart of the argument lies in exhibiting suitable oldforms in \S\ref{oldforms}.\footnote{The reading of \cite{Edi06}*{\S2} was beneficial for the genesis of \S\ref{oldforms}.} 

Ultimately, the sought oldforms come from the analysis of the  degeneracy maps 
\[
\tst \pi_\forg,\, \pi_\quot \colon X_1(n) \ra X_1(\f{n}{2}) \qq \text{over} \q \bF_2 \q \text{for an}  \q n \q \text{with} \q \ord_2(n) = 1
\]
(and their analogues for $X_0(n)$), but at the cost of several complications. Firstly, to exploit the moduli interpretations and to overcome the failure of (S$_2$) of $\Omega^1_{X_1(n)/\bZ_{(2)}}$, we are forced to work with the line bundle $\omega^{\tensor 2}(-\cusps)$ of weight $2$ cusp forms, and hence also with the $\Gamma_1(n)$-level (resp.,~$\Gamma_0(n)$-level) modular stack $\sX_1(n)$ (resp.,~$\sX_0(n)$) in place of its coarse moduli scheme (albeit the difference only manifests itself for $\Gamma_0(n)$). The passage to stacks is facilitated by a certain comparison result overviewed and proved in \Cref{append}. Secondly, several key arguments rest on intersection theory for $\sX_0(n)$ and $\sX_1(n)$, so we crucially use the regularity of these stacks (which may fail for coarse spaces). At multiple places of the overall proof, the moduli interpretations and the analysis of $\sX_0(n)$ and $\sX_1(n)$ presented in \cite{Ces15a} come in handy---although we primarily work over $\bZ_{(2)}$, we cannot ignore the subtleties of the moduli interpretation of $\sX_0(n)_{\bZ_{(2)}}$ at the cusps caused by the fact that $n$ may be divisible by the square of an odd prime.

In the case of \Cref{stevens-main}, the resulting proof is \emph{a posteriori} carried out entirely with schemes because the relevant stacks $\sX_1(n)_{\bZ_{(2)}}$ and $\sX_1(\f{n}{2})_{\bZ_{(2)}}$ identify with their coarse spaces (see \eqref{no-stack}). In contrast, we do not know how to carry out the proof of \Cref{manin-main}~\ref{MM-1}--\ref{MM-15} without resorting to stacks. \Cref{manin-main}~\ref{MM-2} is based on a direct reduction to \Cref{stevens-main}.
\epp

\bpp[Notation] \lab{notation}
The following notation will be in place throughout the paper (see also \S\ref{conv}):
\begin{itemize}
\item
For an open subgroup $H \subset \GL_2(\wh{\bZ})$, we let $\sX_H$ denote the level $H$ modular $\bZ$-stack  defined in \cite{DR73}*{IV.3.3} via normalization (so $\sX_H$ is always Deligne--Mumford and is a scheme for ``small enough'' $H$; see \cite{Ces15a}*{\S4.1} for a review of basic properties of $\sX_H$);

\item
We let $X_H$ denote the coarse moduli space of $\sX_H$, so $X_H$ is the ``usual'' projective modular curve over\footnote{Unlike in the introduction, where $X_0(n)$ and $X_1(n)$ were curves over $\bQ$, all the modular curves and stacks in the rest of the paper are assumed to be over $\bZ$ and we use base change notation $X_0(n)_\bQ$, etc., to denote their $\bQ$-fibers.} $\bZ$ of level $H$ (see \cite{Ces15a}*{6.1--6.3} for a review of basic properties of $X_H$); 

\item
For an $n \in \bZ_{\ge 1}$, we let $\Gamma_0(n) \subset \GL_2(\wh{\bZ})$ (resp.,~$\Gamma_1(n)\subset \GL_2(\wh{\bZ})$) be the preimage of the subgroup $\{ \p{\begin{smallmatrix} * & * \\ 0 & * \end{smallmatrix}} \} \subset \GL_2(\bZ/n\bZ)$ (resp.,~of the subgroup $\{ \p{\begin{smallmatrix} 1 & * \\ 0 & * \end{smallmatrix}} \} \subset \GL_2(\bZ/n\bZ)$);

\item
We write $\sX(1)$ for $\sX_{\GL_2(\wh{\bZ})}$ and sometimes write $\sX_0(n)$ and $\sX_1(n)$ (resp.,~$X_0(n)$ and $X_1(n)$) for $\sX_{\Gamma_0(n)}$ and $\sX_{\Gamma_1(n)}$ (resp.,~for $X_{\Gamma_0(n)}$ and $X_{\Gamma_1(n)}$);

\item
We let $J_0(n) \ce \Pic^0_{X_0(n)_\bQ/\bQ}$ and $J_1(n) \ce \Pic^0_{X_1(n)_\bQ/\bQ}$ be the Jacobian varieties of $X_0(n)_\bQ$ and $X_1(n)_\bQ$, respectively (so $J_0(n)$ and $J_1(n)$ are abelian varieties over $\bQ$);

\item
For a Cohen--Macaulay morphism $\sX \ra S$ (see \S\ref{conv}) of some pure relative dimension from a Deligne--Mumford stack $\sX$ to a scheme $S$, we let $\Omega_{\sX/S}$ (or simply $\Omega$) denote the ``relative dualizing'' quasi-coherent $\sO_\sX$-module discussed in \S\ref{rel-dual}. (We likewise shorten $\Omega^1_{\sX/S}$ to $\Omega^1$.)
\end{itemize}
\epp

\bpp[Conventions] \lab{conv}
A morphism $\sX \ra S$ from a Deligne--Mumford stack (or a scheme) $\sX$ towards a scheme $S$ is \emph{Cohen--Macaulay} if it is flat, locally of finite presentation, and its fibers are Cohen--Macaulay. We write $\sO_{\sX, x}^\sh$ for the strict Henselization of $\sX$ at a geometric point $x$.  

On a modular curve over a subfield of $\bC$, we identify a weight two cusp form with its corresponding K\"{a}hler differential. We use the $j$-invariant to identify $X_{\GL_2(\wh{\bZ})}$ with $\bP^1_\bZ$ (see \cite{DR73}*{VI.1.1 and VI.1.3}). We use `new' and `optimal' in the sense of the beginning of the introduction.

For a proper smooth geometrically connected curve $X$ over a field $k$, we make the identification 
\be \lab{id-omg}
H^0(X, \Omega^1) \cong H^0(\Pic^0_{X/k}, \Omega^1)
\ee
supplied by the combination of Grothendieck--Serre duality and the deformation-theoretic identification $H^1(X, \sO_{X}) \cong \Lie(\Pic^0_{X/k})$, and, whenever we choose an $x_0 \in X(k)$, we freely use the alternative description of the identification \eqref{id-omg} as pullback of K\"{a}hler differentials along the ``$x \mapsto \sO(x) \tensor \sO(x_0)\i$'' closed immersion $X \hra \Pic^0_{X/k}$ (see \cite{Con00}*{Thm.~B.4.1}).

An element of a torsion free module over a Dedekind domain is \emph{primitive} if the quotient by the submodule that it generates is torsion free. For a prime $p$, we let $\ord_p$ denote the $p$-adic valuation with $\ord_p(p) = 1$ and let $(-)_{(p)}$ denote localization at $p$. For an $n \in \bZ_{\ge 1}$, we set $\mu_n \ce \Ker(\bG_m \xra{n} \bG_m)$, let $\bQ(\zeta_n)$ denote the $n\th$ cyclotomic field, let $\bZ[\zeta_n]$ denote its ring of integers, and let $\bZ[\zeta_n]^+$ denote the ring of integers of the maximal totally real subfield of $\bQ(\zeta_n)$. A dual abelian variety, a Cartier dual commutative finite locally free group scheme, or a dual homomorphism is denoted by $(-)^\vee$.
\epp

\subsection*{Acknowledgements}
I thank Brian Conrad for several extended discussions. I thank Kevin Buzzard, Frank Calegari, John Cremona, Bas Edixhoven, Haruzo Hida, Erick Knight, Martin Olsson, Michel Raynaud, Ken Ribet, Sug Woo Shin, Yunqing Tang, Xinyi Yuan, and David Zureick-Brown for helpful conversations or correspondence. I thank the referee for helpful comments. I thank the MathOverflow community---the reading of several anonymous discussions was useful while working on this paper. I thank the Miller Institute for Basic Research in Science at the University of California Berkeley for its support.


\subfile{reduction}


\subfile{oldforms}


\subfile{proofs}






\appendix

\subfile{app-integral}

\end{document}

%% file: reduction.tex

\section{A reduction to a problem about congruences between modular forms} \lab{reduction}

Our approach to the conjectures of Manin and Stevens rests on \Cref{congJ-deg}~\ref{CD-c}, which relates them to a comparison between the modular degree and a certain congruence number. Our first task is to define the latter in \S\ref{ref-cong} after introducing the relevant setup in \S\ref{setup}.

\bpp[The setup] \lab{setup} Throughout \S\ref{reduction} we supplement the notation of \S\ref{notation} with the following:
\begin{itemize}
\item
We let $\Gamma$ denote either $\Gamma_0(n)$ or $\Gamma_1(n)$ for a fixed $n \in \bZ_{\ge 1}$;

\item
We let $X$ denote $X_\Gamma$, i.e., either $X = X_0(n)$ or $X = X_1(n)$;

\item
We let $J$ denote the Jacobian $\Pic^0_{X_\bQ/\bQ}$, i.e., either $J = J_0(n)$ or $J = J_1(n)$;

\item
We let $\pi\colon J \surjects E$ be a new elliptic optimal quotient (so $E$ is an elliptic curve over $\bQ$ and $\Ker \pi$ is an abelian variety);

\item
We let $\pi \colon \cJ \ra \cE$ denote the extension to N\'{e}ron models over $\bZ$;

\item
We let $f \in H^0(X_\bQ, \Omega^1) \overset{\eqref{id-omg}}{\cong} H^0(J, \Omega^1)$ denote the normalized new eigenform corresponding~to $\pi$ (`normalized' means that the $q$-expansion $\p{\sum_{n \ge 1} a_n q^n} \f{dq}{q}$ of $f$ at the cusp\footnote{We form $q$-expansions after identifying $(X_\Gamma)_\bC$ with the quotient of the upper half-plane, so over $\bC$. When $\Gamma = \Gamma_1(n)$, the cusp $\infty$ does not descend to a $\bQ$-point of $X$, and the $q$-expansion of $f$ need not have rational coefficients.} ``$\infty$'' has $a_1 = 1$);

\item
We let $c_\pi \in \bQ^\times$ denote the Manin(--Stevens) constant of $E$, i.e., $\pi^*(\omega_\cE) = c_\pi \cdot f$ in $H^0(J, \Omega^1)$, where $\omega_\cE \in H^0(E, \Omega^1)$ is a generator of $H^0(\cE, \Omega^1)$ (so $c_\pi$ is only well-defined up to $\pm 1$).

\end{itemize}
One knows that $c_\pi \in \bZ$ (see \cite{Edi91}*{Prop.~2} and \cite{Ste89}*{Thm.~1.6}), and the conjectures of Manin and Stevens predict that $c_\pi = \pm 1$. 
\epp

\bpp[Congruence with respect to the lattice $H^0(\cJ, \Omega^1)$] \lab{ref-cong} The relevant ``congruence module'' is
\be \lab{cong-mod}
\tst \f{H^0(\cJ, \Omega^1)}{H^0(\cJ, \Omega^1) \cap (\bQ \cdot f) + H^0(\cJ, \Omega^1) \cap (\bQ \cdot f)^\perp},
\ee
where the orthogonal complement is taken in $H^0(X_\bQ, \Omega^1)$ with respect to the Petersson inner product. This $\bZ$-module is finite and cyclic (because $\f{H^0(\cJ, \Omega^1)}{H^0(\cJ, \Omega^1) \cap (\bQ \cdot f)^\perp} \simeq \bZ$), and we denote its order by 
\[
\tst \mathrm{cong}_{f,\, \cJ} \ce \#\p{\f{H^0(\cJ, \Omega^1)}{H^0(\cJ, \Omega^1) \cap (\bQ \cdot f) + H^0(\cJ, \Omega^1) \cap (\bQ \cdot f)^\perp}}.
\]
\epp

\bpropt \lab{congJ-deg}
Let $\phi\colon X_\bQ \surjects E$ denote the composition of $\pi\colon J \surjects E$ with the immersion $i_P\colon X_\bQ \hra J$ obtained by choosing a base point $P \in X(\bQ)$ \up{for instance, a rational cusp}.
\benum
\item \lab{CD-a}
The composition $\pi \circ \pi^\vee \colon E \ra J \ra E$ is multiplication by $\deg \phi$ \up{which is independent of $P$}.

\item \lab{CD-b}
With the notation of \uS\uref{ref-cong},
\[
\q \mathrm{cong}_{f,\, \cJ} \mid \deg \phi.
\]

\item \lab{CD-c}
If $p$ is a prime such that $f \in H^0(\cJ_{\bZ_{(p)}}, \Omega^1)$, then 
\[
\tst \q  \ord_p(c_\pi) \le \ord_p(\f{\deg \phi}{\mathrm{cong}_{f,\, \cJ}}).
\]
\eenum
\epropt

\bpf \hfill
\benum
\item
We compute the effect of $\pi \circ \pi^\vee$ on a variable point $Q \in E(\ov{\bQ})$. The canonical principal polarization of $E$ sends $Q$ to $\sO_{E_{\ov{\bQ}}}([-Q] - [0])$, which $\Pic^0(\phi) = \Pic^0(i_P) \circ \pi^\vee$ then sends to $\sO_{X_{\ov{\bQ}}}([\phi\i(-Q)] - [\phi\i(0)])$. Thus, since $\Pic^0(i_P)$ is the negative of the inverse of the canonical principal polarization of $J$ (see, for instance, \cite{Mil86}*{6.9}), the overall effect of $\pi \circ \pi^\vee$ is to send $Q$ to the negative of the sum under the group law of $E_{\ov{\bQ}}$ of the $\phi$-image of the divisor $[\phi\i(-Q)] - [\phi\i(0)]$ on $X_{\ov{\bQ}}$, i.e.,~to $\deg \phi \cdot Q$.

\item
Due to the optimality of $\pi$, the dual $\pi^\vee\colon E \ra J$ is a closed immersion. Since $\pi^\vee$ is Hecke equivariant (see \cite{ARS12}*{pp.~24--25}), it induces the injection
\be \lab{pull-inj}
\tst \qq \f{H^0(\cJ, \Omega^1)}{H^0(\cJ, \Omega^1) \cap (\bQ \cdot f)^\perp} \hra H^0(\cE, \Omega^1).
\ee
Due to the Hecke equivariance of $\pi$, we have $\pi^*(H^0(\cE, \Omega^1)) \subset H^0(\cJ, \Omega^1) \cap (\bQ \cdot f)$. Moreover, the $\bZ$-line $H^0(\cJ, \Omega^1) \cap (\bQ \cdot f)$ maps injectively into the source of \eqref{pull-inj}, so from \ref{CD-a} and \eqref{pull-inj} we get the injection
\be \lab{cong-inj}
\tst \qq \f{H^0(\cJ, \Omega^1)}{\pi^*(H^0(\cE, \Omega^1)) + H^0(\cJ, \Omega^1) \cap (\bQ \cdot f)^\perp} \hra \f{H^0(\cE, \Omega^1)}{(\deg \phi) \cdot H^0(\cE, \Omega^1)}
\ee
that exhibits the ``congruence module'' of \eqref{cong-mod} as a subquotient of $\bZ/(\deg \phi)\bZ$. It remains to observe that the order of every subquotient of $\bZ/(\deg \phi)\bZ$ divides $\deg \phi$.

\item
By quantifying at a prime $p$ the extent to which the inclusion \eqref{cong-inj} fails to be an isomorphism between $\f{H^0(\cE, \Omega^1)}{(\deg \phi) \cdot H^0(\cE, \Omega^1)}$ and the ``congruence module'' of \eqref{cong-mod}, we arrive at the equality
\[
\tst \qq \q \ord_p\p{\f{\deg \phi}{\mathrm{cong}_{f,\, \cJ}}} = \ord_p\p{\#\p{\f{H^0(\cE, \Omega^1)}{\im((\pi^\vee)^*\colon H^0(\cJ, \Omega^1) \ra H^0(\cE, \Omega^1))}}} + \ord_p\p{\#\p{\f{H^0(\cJ, \Omega^1) \cap (\bQ \cdot f)}{\pi^*(H^0(\cE, \Omega^1))}}}.
\]
Since $\bZ_{(p)} \cdot f \subset H^0(\cJ_{\bZ_{(p)}}, \Omega^1) \cap (\bQ \cdot f)$ and $\pi^*(H^0(\cE, \Omega^1)) = \bZ \cdot c_\pi f$ with $c_\pi \in \bZ$, the last summand is at least $\ord_p(c_\pi)$, and the sought inequality follows.
 \qedhere
\eenum
\epf

Applying \Cref{congJ-deg}~\ref{CD-c} to study the conjectures of Manin and Stevens at a prime $p$ essentially amounts to establishing the $p$-part of the divisibility converse to the divisibility $\mathrm{cong}_{f,\, \cJ} \mid \deg \phi$ supplied by \Cref{congJ-deg}~\ref{CD-b}. A result of Ribet, \cite{ARS12}*{Thm.~3.6~(a)} (see also \cite{AU96}*{Lem.~3.2} and \cite{CK04}*{Thm.~1.1} for other expositions in the case when $\Gamma = \Gamma_0(n)$), supplies the sought converse divisibility, but with the caveat that the congruences be considered with respect to another lattice $S_2(\Gamma, \bZ) \subset H^0(X_\bQ, \Omega^1)$ in place of $H^0(\cJ, \Omega^1)$. Therefore, our task is to relate the two types of congruences. For this, we work under the assumption that $\ord_p(n) \le 1$ and focus on the key case $p = 2$ (although, as we point out along the way, for $\Gamma_0(n)$ most of the arguments also work for odd $p$). In this setting, we relate the two types of congruences in the proofs of \Cref{odd-n,main-reduction}. 

In the focal case $p = 2$ with $\ord_2(n) \le 1$, we begin with the simpler possibility $\ord_2(n) = 0$.

\bpp[The structure of $X_{\bZ_{(p)}}$ when $\ord_p(n) = 0$] \lab{str-odd-n}
If $p \nmid n$, then the ``level'' of $\Gamma$ is prime to $p$, so the proper $\bZ_{(p)}$-curve $X_{\bZ_{(p)}}$ is  smooth (see \cite{DR73}*{VI.6.7}, possibly also \cite{Ces15a}*{6.4~(a)}). Moreover, its geometric fibers are irreducible by \cite{DR73}*{IV.5.6}. Due to these properties, $\Pic^0_{X_{\bZ_{(p)}}/\bZ_{(p)}}$ is an abelian $\bZ_{(p)}$-scheme (see \cite{BLR90}*{9.4/4}), and hence identifies with $\cJ_{{\bZ_{(p)}}}$. In particular,
\be \lab{int-comp}
H^0(X_{\bZ_{(p)}}, \Omega^1) = H^0(\cJ_{{\bZ_{(p)}}}, \Omega^1) \qq \text{inside} \qq H^0(X_{\bQ}, \Omega^1).
\ee
\epp

The method of proof of the following theorem in essence amounts to the method used in \cite{AU96}*{proof of Thm.~A} in the setting of $\Gamma_0(n)$. At least when $\Gamma = \Gamma_0(n)$, the method is not specific to $p = 2$.

\bthmt \lab{odd-n}
If $\ord_2(n) = 0$, then, in the setting of \uS\uS\uref{setup}--\uref{ref-cong},
\[
\ord_2(\mathrm{cong}_{f,\, \cJ}) = \ord_2(\deg \phi) \qq \text{and} \qq \ord_2(c_\pi) = 0.
\]
\ethmt

\bpf
For a subring $R \subset \bC$, let $S_2(\Gamma, R)$ denote the $R$-module of those weight $2$ cusp forms of level $\Gamma$ whose Fourier expansion at the cusp ``$\infty$'' has coefficients in $R$. As described in \cite{DI95}*{\S12.3},
\be\lab{Shi-input}
S_2(\Gamma, R) \cong S_2(\Gamma, \bZ) \tensor_\bZ R
\ee
(both sides of this identification are defined in terms of $(X_\Gamma)_\bC$, so for arguing it in the $\Gamma = \Gamma_1(n)$ case one is free to use the ``$\mu$-model'' of the modular curve to descend the cusp ``$\infty$'' to a $\bZ$-point).

If $\Gamma = \Gamma_0(n)$, then \eqref{int-comp} and \cite{Edi06}*{2.2 and 2.5} ensure that 
\[
H^0(\cJ_{\bZ_{(2)}}, \Omega^1) = S_2(\Gamma, \bZ_{(2)}) \qq \text{inside}\qq H^0(X_\bQ, \Omega^1) = S_2(\Gamma, \bQ).
\]
Thus, in the $\Gamma = \Gamma_0(n)$ case,
\be \lab{cong-comp-1}
\tst \p{\f{S_2(\Gamma, \bZ)}{S_2(\Gamma, \bZ) \cap \bQ \cdot f + S_2(\Gamma, \bZ) \cap (\bQ \cdot f)^\perp}} \tensor_\bZ \bZ_{(2)} \cong \p{\f{H^0(\cJ, \Omega^1)}{H^0(\cJ, \Omega^1) \cap \bQ \cdot f + H^0(\cJ, \Omega^1) \cap (\bQ \cdot f)^\perp}} \tensor_\bZ \bZ_{(2)}.
\ee
By \cite{ARS12}*{Thm.~3.6~(a)}, the order of the left side of \eqref{cong-comp-1} is divisible by $2^{\ord_2(\deg \phi)}$, so $\ord_2(\deg \phi) \le \ord_2(\mathrm{cong}_{f,\, \cJ})$. Due to the converse inequality supplied by \Cref{congJ-deg}~\ref{CD-b}, equality must hold. \Cref{congJ-deg}~\ref{CD-c} then settles the $\Gamma = \Gamma_0(n)$ case because $c_\pi \in \bZ$ and the equality $H^0(\cJ_{\bZ_{(2)}}, \Omega^1) = S_2(\Gamma, \bZ_{(2)})$ also provides the containment $f \in H^0(\cJ_{\bZ_{(2)}}, \Omega^1)$.

For the remainder of the proof we assume that $\Gamma = \Gamma_1(n)$. One special feature of this case is that $X_{\bZ_{(2)}} = \sX_{\bZ_{(2)}}$ due to the triviality of the automorphism functors of the geometric points of $\sX_{\bZ_{(2)}}$ forced by the inequality $n \ge 5$ resulting from the existence of $f$ (see \cite{Ces15a}*{4.1.4, 4.4.4~(c)} and \cite{KM85}*{2.7.4}). By pulling back $f$ along the $\bZ_{(2)}$-base change of the forgetful map $X_1(n) \ra X_0(n)$, we see with the help of \eqref{int-comp} that the containment $f \in H^0(\cJ_{\bZ_{(2)}}, \Omega^1)$ continues to hold.

The cusp ``$\infty$'' arises from a $\bZ[\zeta_n]$-point (even a $\bZ[\zeta_n]^+$-point) of $X$ via an embedding $\bQ(\zeta_n) \subset \bC$ whose choice we fix, and the completion of $X_{\bZ[\zeta_n]_{(2)}}$ along the resulting $\bZ[\zeta_n]_{(2)}$-point is isomorphic to $\bZ[\zeta_n]_{(2)}\llb q \rrb$ and is described by a Tate curve (combine \cite{DR73}*{VII.2.1} and \cite{KM85}*{1.12.9}; see also \cite{Con07a}*{4.3.7}). 
Therefore, \cite{Edi06}*{proof of 2.2 and top of p.~6} provide the identification 
\[
H^0(X_{\bQ(\zeta_n)}, \Omega^1) = S_2(\Gamma, \bQ(\zeta_n)) \qq \text{under which} \qq H^0(X_{\bZ[\zeta_n]_{(2)}}, \Omega^1) = S_2(\Gamma, \bZ[\zeta_n]_{(2)}).
\]
Therefore, with the help of \eqref{int-comp} and \eqref{Shi-input} we obtain the following analogue of \eqref{cong-comp-1}:
\be\lab{cong-comp-2}
\tst \p{\f{S_2(\Gamma, \bZ)}{S_2(\Gamma, \bZ) \cap \bQ \cdot f + S_2(\Gamma, \bZ) \cap (\bQ \cdot f)^\perp}} \tensor_\bZ \bZ[\zeta_n]_{(2)} \cong \p{\f{H^0(\cJ, \Omega^1)}{H^0(\cJ, \Omega^1) \cap \bQ \cdot f + H^0(\cJ, \Omega^1) \cap (\bQ \cdot f)^\perp}} \tensor_\bZ \bZ[\zeta_n]_{(2)}.
\ee
By \cite{ARS12}*{Thm.~3.6 (a)}, the exponent of the left side of \eqref{cong-comp-2} is divisible by $\ord_2(\deg \phi)$. Since the exponent of the right side equals $\ord_2(\mathrm{cong}_{f,\, \cJ})$, the resulting inequality 
\[
\ord_2(\deg \phi) \le \ord_2(\mathrm{cong}_{f,\, \cJ})
\]
combines with \Cref{congJ-deg}~\ref{CD-b}--\ref{CD-c} to conclude the proof as in the case when $\Gamma = \Gamma_0(n)$.
\epf

We turn our attention to the more complex possibility $\ord_2(n) = 1$.

\bpp[The structure of $X_{\bZ_{(2)}}$ when $\ord_2(n) = 1$] \lab{spec-fib}
The $\bZ_{(2)}$-curve $X_{\bZ_{(2)}}$ is always normal, and proper and flat over $\bZ_{(2)}$ (these are general properties of modular curves, see, for instance, \cite{Ces15a}*{6.1}). Moreover, the $\bZ_{(2)}$-fibers of $X_{\bZ_{(2)}}$ are geometrically connected by \cite{DR73}*{IV.5.5}. However, if $\ord_2(n) = 1$, as we assume from now on, then the $\bF_2$-fiber is singular. Nevertheless, if we set $\Gamma' \ce \Gamma_0(\f{n}{2})$ when $\Gamma = \Gamma_0(n)$ and $\Gamma' \ce \Gamma_1(\f{n}{2})$ when $\Gamma = \Gamma_1(n)$, then we have
\be \lab{no-2}
\Gamma = \Gamma_0(2) \cap \Gamma'.
\ee
Therefore, by \cite{DR73}*{VI.6.9}, $X_{\bF_2}$ is semistable and has two irreducible components that meet precisely at the supersingular points. Both components are isomorphic to the proper, smooth, geometrically connected $\bF_2$-curve $(X_{\Gamma'})_{\bF_2}$, so the semistable $\bZ_{(2)}$-curve $X_{\bZ_{(2)}}$ is smooth away from the supersingular points on its special fiber. 

The locus of $X_{\bF_2}$ that corresponds to ordinary elliptic curves is a disjoint union of two affine connected curves: the open whose geometric points correspond to $\Gamma$-level structures with a connected $\Gamma_0(2)$-part, and the open for which this $\Gamma_0(2)$-part is \'{e}tale. We let 
\[
X_{\bF_2}^\mu \qq \text{(resp., } X_{\bF_2}^\et)
\]
denote the irreducible component of $X_{\bF_2}$ that contains the former (resp.,~the latter) open, and we define the $\bZ_{(2)}$-smooth open $U^\mu \subset X_{\bZ_{(2)}}$ by
\[
U^\mu \ce X_{\bZ_{(2)}} \setminus X_{\bF_2}^\et.
\]

The existence of $f$ ensures that $\f{n}{2} \ge 5$, so, as in the proof of \Cref{odd-n},
\[
X_{\bZ_{(2)}} = \sX_{\bZ_{(2)}} \qq \text{if} \qq \Gamma = \Gamma_1(n).
\]
In particular, $X_{\bZ_{(2)}}$ is regular when $\Gamma = \Gamma_1(n)$ (but need not be regular when $\Gamma = \Gamma_0(n)$), see \cite{Ces15a}*{4.4.4}. The semistability of $X_{\bZ_{(2)}}$ supplies the identification 
\be \lab{pic-jic}
\Pic^0_{X_{\bZ_{(2)}}/\bZ_{(2)}} \cong \cJ_{\bZ_{(2)}}^0
\ee
as in \cite{BLR90}*{9.7/2} and ensures that the relative dualizing sheaf $\Omega$ is a line bundle on $X_{\bZ_{(2)}}$. In particular, \eqref{pic-jic} and Grothendieck duality as in \cite{Con00}*{Cor.~5.1.3} supply the analogue of \eqref{int-comp}:
\be\lab{int-comp-2}
H^0(X_{\bZ_{(2)}}, \Omega) = H^0(\cJ_{{\bZ_{(2)}}}, \Omega^1) \qq \text{inside} \qq H^0(X_{\bQ}, \Omega^1).
\ee
Although $U^\mu$ is not $\bZ_{(2)}$-proper, $H^0(U^\mu, \Omega^1)$ is a finite free $\bZ_{(2)}$-module that contains $H^0(X_{\bZ_{(2)}}, \Omega)$ and identifies with a $\bZ_{(2)}$-lattice inside $H^0(X_{\bQ}, \Omega^1)$, see \cite{BDP16}*{Prop.~B.2.1.1}.
\epp

When combined with \eqref{int-comp-2}, the following lemma will aid our analysis of the congruence module \eqref{cong-mod} in the case when $\ord_2(n) = 1$. The involution trick used in its proof may be traced back at least to \cite{Maz78}*{proof of Prop.~3.1}.

\blemt \lab{hop-f}
If $\ord_2(n) = 1$, then, in the setting of \uS\uref{setup} and \uS\uref{spec-fib},
\be \lab{hop-f-eq}
H^0(U^\mu, \Omega^1) \cap (\bQ \cdot f) = H^0(X_{\bZ_{(2)}}, \Omega) \cap (\bQ \cdot f) \qq \text{inside} \qq H^0(X_\bQ, \Omega^1),
\ee
and $f \in H^0(X_{\bZ_{(2)}}, \Omega) \overset{\eqref{int-comp-2}}{=} H^0(\cJ_{{\bZ_{(2)}}}, \Omega^1)$.
\elemt

\bpf
Since $\Omega$ is a line bundle, the normality of $X_{\bZ_{(2)}}$ ensures that a $g \in H^0(X_\bQ, \Omega^1)$ extends to $H^0(U^\mu, \Omega^1)$ if and only if $g$ extends to the stalk of $\Omega^1_{X_{\bZ_{(2)}}/\bZ_{(2)}}$ at the generic point of $U^\mu_{\bF_2}$, in which case $g$ extends further to $H^0(X_{\bZ_{(2)}}, \Omega)$ if and only if it extends to the stalk of $\Omega^1_{X_{\bZ_{(2)}}/\bZ_{(2)}}$ at the other generic point of $X_{\bF_{2}}$. Due to \eqref{no-2}, the Atkin--Lehner involution $w_2$ makes sense on the elliptic curve locus of $X_{\bZ_{(2)}}$. Moreover, it extends to $X_\bQ$ and interchanges the two stalks considered above. The equality \eqref{hop-f-eq} follows because the effect of $w_2$ on the newform $f$ is scaling by $\pm 1$. 

Due to \eqref{hop-f-eq}, for the rest it suffices to note that $f \in H^0(U^\mu, \Omega^1)$ in the case when $\Gamma = \Gamma_0(n)$ (see \cite{Edi06}*{2.5}), and hence also in general thanks to the forgetful map $X_{\Gamma_1(n)} \ra X_{\Gamma_0(n)}$.
\epf

\bremst 
\remit \lab{also-odd-p}
In the case $\Gamma = \Gamma_0(n)$, the discussion of \S\ref{spec-fib} and the proof of \Cref{hop-f} are not specific to the prime $p = 2$. In particular, they show that if $\Gamma = \Gamma_0(n)$, then for any prime $p$ with $\ord_p(n) \le 1$ we have 
\[
\qq f \in H^0(X_{\bZ_{(p)}}, \Omega) = H^0(\cJ_{\bZ_{(p)}}, \Omega^1).
\]
Continuing to assume that $p$ is a prime with $\ord_p(n) \le 1$, we claim that this implies that
\[
\qq f \in H^0(\cJ_{\bZ_{(p)}}, \Omega^1) \qq  \text{also when $\Gamma = \Gamma_1(n)$.}
\]
For this, we first fix an $x \in X_1(n)(\bQ)$ and consider the resulting immersions $X_1(n) \hra J_1(n)$ and $X_0(n) \hra J_0(n)$, which, by Albanese functoriality, induce a compatible homomorphism $J_1(n) \ra J_0(n)$. To then see the claim, it remains to pullback $f$ along the resulting map of N\'{e}ron models over $\bZ_{(p)}$ and to use the alternative description of \eqref{id-omg} mentioned in \S\ref{conv}.

\remit \lab{degree-divides}
In the $\ord_2(n) = 1$ case, \eqref{int-comp-2} and \Cref{hop-f} guarantee that $f \in H^0(\cJ_{\bZ_{(2)}}, \Omega^1)$, so \Cref{congJ-deg}~\ref{CD-c} supplies the inequality $\ord_2(c_\pi) \le \ord_2\p{\f{\deg \phi}{\mathrm{cong}_{f,\, \cJ}}}$. In particular, 
\[
\text{if} \q \ord_2(n) = 1 \q \text{and}  \q 2\nmid \deg \phi, \q \text{then} \q \ord_2(c_\pi) = 0,
\]
which for $\Gamma = \Gamma_0(n)$ recovers a result of Agashe--Ribet--Stein, \cite{ARS06}*{Thm.~2.7} (see also \cite{ARS06}*{Thm.~3.11} for a generalization to optimal newform quotients of arbitrary dimension).

\eremst

In the $\ord_2(n) = 1$ case, the main result of this section is the following outgrowth of \Cref{congJ-deg}~\ref{CD-c}.

\bthmt \lab{main-reduction}
If $\ord_2(n) = 1$, then, in the setting of \uS\uref{setup} and \uS\uref{spec-fib}, the group
\[
\tst \f{H^0(U^\mu, \Omega^1)}{H^0(X_{\bZ_{(2)}}, \Omega) + H^0(U^\mu, \Omega^1) \cap (\bQ \cdot f)^\perp}
\]
is finite and cyclic and the $2$-adic valuation of its order bounds the $2$-adic valuation of $c_\pi$\ucolon
\[
\tst \ord_2(c_\pi) \le \ord_2 \p{\#\p{\f{H^0(U^\mu, \Omega^1)}{H^0(X_{\bZ_{(2)}}, \Omega) + H^0(U^\mu, \Omega^1) \cap (\bQ \cdot f)^\perp}}}.
\]
\ethmt

\emph{Proof.}
Due to \Cref{hop-f}, the pullback map
\be \lab{key-comp}
\tst \f{H^0(X_{\bZ_{(2)}}, \Omega)}{H^0(X_{\bZ_{(2)}}, \Omega) \cap (\bQ \cdot f) + H^0(X_{\bZ_{(2)}}, \Omega) \cap (\bQ \cdot f)^\perp} \hra \f{H^0(U^\mu, \Omega^1)}{H^0(U^\mu, \Omega^1) \cap (\bQ \cdot f) + H^0(U^\mu, \Omega^1) \cap (\bQ \cdot f)^\perp}
\ee
is injective and its cokernel is the group in question, which therefore inherits finiteness and cyclicity from the target (compare with the discussion of finiteness and cyclicity in \S\ref{ref-cong}). 

The sought inequality follows by combining \Cref{congJ-deg}~\ref{CD-c}, \Cref{hop-f}, and the following claims.

\bcl \lab{claim-source}
The $2$-adic valuation of the order of the source of \eqref{key-comp} is $\ord_2(\mathrm{cong}_{f,\, \cJ})$.
\ecl

\bcl \lab{claim-target}
The $2$-adic valuation of the order of the target of \eqref{key-comp} is at least $\ord_2(\deg \phi)$.
\ecl

\bpf[Proof of Claim \uref{claim-source}]
It suffices to use the identification \eqref{int-comp-2}.
\epf

\emph{Proof of Claim \uref{claim-target}.}
We use the same notation $S_2(\Gamma, R)$ as in the proof of \Cref{odd-n}. Like there, for every subring $R \subset \bC$ we have the identification
\[
S_2(\Gamma, R) \cong S_2(\Gamma, \bZ) \tensor_\bZ R
\]
discussed in \cite{DI95}*{\S12.3}.

If $\Gamma = \Gamma_0(n)$, then $H^0(U^\mu, \Omega^1) = S_2(\Gamma, \bZ_{(2)})$ (see \cite{Edi06}*{2.5}), so \cite{ARS12}*{Thm.~3.6 (a)} shows that the $2$-primary factor of $\deg \phi$ divides the order of the target of \eqref{key-comp}.

If $\Gamma = \Gamma_1(n)$, then we carry out the same argument after base change to $\bZ[\zeta_n]_{(2)}$. Namely, after fixing an embedding $\bQ(\zeta_n) \subset \bC$, we arrive at the identification
\[
H^0(U^\mu_{\bZ[\zeta_n]_{(2)}}, \Omega^1) = S_2(\Gamma, \bZ[\zeta_n]_{(2)}),
\]
analogously to the proof of \Cref{odd-n}. In particular, the base change to $\bZ[\zeta_n]_{(2)}$ of the target of \eqref{key-comp} identifies with 
\[
\tst \p{\f{S_2(\Gamma, \bZ)}{S_2(\Gamma, \bZ) \cap (\bQ \cdot f) + S_2(\Gamma, \bZ) \cap (\bQ \cdot f)^\perp}} \tensor_\bZ \bZ[\zeta_n]_{(2)},
\]
so it remains to observe that the exponent of the latter is divisible by the $2$-primary factor of $\deg \phi$ due to \cite{ARS12}*{Thm.~3.6 (a)}.
\QEDD

\bremt \lab{more-odd-p}
In the case $\Gamma = \Gamma_0(n)$, neither the statement nor the proof of \Cref{main-reduction} is specific to the prime $p = 2$, as one sees with the help of  Remark \ref{also-odd-p}.
\eremt

%% file: oldforms.tex

\section{Using oldforms to offset the difference between integral structures} \lab{oldforms}

According to \Cref{main-reduction}, to settle the $\ord_2(n) = 1$ case of the $2$-primary part of \Cref{stevens-conj,manin-conj}, it suffices to show that $H^0(U^\mu, \Omega^1)/H^0(X_{\bZ_{(2)}}, \Omega)$ consists of images of elements of $H^0(U^\mu, \Omega^1) \cap (\bQ \cdot f)^\perp$. Since $f$ is a newform, the latter space contains the oldforms in $H^0(U^\mu, \Omega^1)$, and our strategy is to show that under suitable assumptions these oldforms sweep out the entire $H^0(U^\mu, \Omega^1)/H^0(X_{\bZ_{(2)}}, \Omega)$. The merit of this approach is that the sought statement no longer involves newforms, but instead is a generality about integral structures on the $\bQ$-vector space of weight $2$ cusp forms. We therefore forget about $f$ and pursue this generality with the following setup.

\bpp[The setup] \lab{setup-2}
Throughout this section we fix an $n\in \bZ_{\ge 1}$ with $\ord_2(n) = 1$ and
\begin{itemize}
\item
We let $\Gamma$ and $\Gamma'$ denote either $\Gamma_0(n)$ and $\Gamma_0(\f{n}{2})$, respectively, or $\Gamma_1(n)$ and $\Gamma_1(\f{n}{2})$, respectively;

\item
We let $\sX$ (resp.,~$\sX'$) denote $\sX_{\Gamma}$ (resp.,~$\sX_{\Gamma'}$), so that, e.g.,~$\sX'$ is either $\sX_0(\f{n}{2})$ or $\sX_1(\f{n}{2})$;

\item
We let $X$ and $X'$ denote the coarse moduli schemes of $\sX$ and $\sX'$, i.e., $X = X_\Gamma$ and $X' = X_{\Gamma'}$;

\item
We set $U^\mu \ce X_{\bZ_{(2)}} \setminus X_{\bF_2}^\et$ and $U^\et \ce X_{\bZ_{(2)}} \setminus X_{\bF_2}^\mu$ (see \S\ref{spec-fib} for the definition of $X_{\bF_2}^\mu$ and $X_{\bF_2}^\et$);

\item
We let $\sU^\mu \subset \sX_{\bZ_{(2)}}$ (resp.,~$\sU^\et \subset \sX_{\bZ_{(2)}}$) be the preimage of $U^\mu$ (resp.,~of $U^\et$).
\end{itemize}

The algebraic stacks $\sX$ and $\sX'$ are regular and have moduli interpretations in terms of generalized elliptic curves equipped with additional data, see \cite{Ces15a}*{\S4.4, esp.~4.4.4, and Ch.~5, esp.~\S\S5.9--5.10 and 5.13~(a)}. These moduli interpretations and \cite{KM85}*{2.7.4} show that
\be \lab{no-stack}
\tst \sX_{\bZ_{(2)}} = X_{\bZ_{(2)}} \qq \text{if} \q \Gamma = \Gamma_1(n) \q \text{with} \q  \f{n}{2} \ge 5
\ee
(see \cite{Ces15a}*{4.1.4}). Even though we do not rely on \eqref{no-stack} in what follows, its significance is that the overall proof of \Cref{stevens-main} is actually carried out entirely in the realm of schemes.
\epp

Using the moduli interpretations of $\sX$ and $\sX'$, we seek to expose degeneracy maps
\[
\tst \pi_\forg,\, \pi_\quot \colon \sX \rightrightarrows \sX'
\]
whose base changes to $\bZ_{(2)}$ will be instrumental for constructing enough oldforms in $H^0(U^\mu, \Omega^1)$. The construction of $\pi_\forg$ and $\pi_\quot$ is not specific to the prime $2$, so we present it in the setting of any prime $p$ and any $N \in \bZ_{\ge 1}$ with $\ord_p(N) = 1$.

\bpp[The maps $\pi_\forg$ and $\pi_\quot$ in the $\Gamma_1(N)$ case] \lab{forg-quot-1}
The stack $\sX_1(N)$ (resp.,~$\sX_1(\f{N}{p})$) parametrizes generalized elliptic curves $E \ra S$ equipped with an ample Drinfeld $\bZ/N\bZ$-structure (resp.,~$\bZ/\f{N}{p}\bZ$-structure). The map 
\[
\tst \pi_\forg\colon \sX_1(N) \ra \sX_1(\f{N}{p})
\]
is defined by forgetting the $p$-primary part of the $\bZ/N\bZ$-structure and contracting $E$ with respect to the $\f{N}{p}$-primary part, whereas the map 
\[
\tst \pi_\quot \colon \sX_1(N) \ra \sX_1(\f{N}{p})
\]
is defined by replacing $E$ by the quotient by the subgroup generated by this $p$-primary part and equipping the quotient with the induced ample Drinfeld $\bZ/\f{N}{p}\bZ$-structure (without restricting to the elliptic curve locus the quotient (here and below) is in the sense of \cite{Ces15a}*{2.2.4 and 2.2.6} and carries an induced $\bZ/\f{N}{p}\bZ$-structure by \cite{Ces15a}*{4.2.9 (e)}).

The maps $\pi_\forg$ and $\pi_\quot$ are representable because, due to the representability of the forgetful contraction $\sX_1(N) \ra \sX(1)$, they do not identify any distinct automorphisms of any geometric point of $\sX_1(N)$ (see \cite{Ces15a}*{3.2.2~(b) and proof of Lemma 4.7.1}). They inherit properness from $\sX_1(N) \ra \Spec \bZ$ and, due to their moduli interpretation, quasi-finiteness from $\sX_1(N) \ra \sX(1)$ (the quasi-finiteness of the latter ensures that over a fixed algebraically closed field there are only finitely many isomorphism classes of degenerate generalized elliptic curves equipped with an ample Drinfeld $\bZ/N\bZ$-structure). Therefore, $\pi_\forg$ and $\pi_\quot$ are finite (see \cite{LMB00}*{A.2} and \cite{EGAIV4}*{18.12.4}), and hence even locally free due to the miracle flatness theorem (see \cite{EGAIV2}*{6.1.5}). 

The same argument will show that the maps $\pi_\forg$ and $\pi_\quot$ are also representable and finite locally free in the $\Gamma_0(N)$ case because $\pi_\forg$ (resp.,~$\pi_\quot$) will still send underlying generalized elliptic curves to their suitable contractions (resp.,~quotients). 
\epp

\bpp[The maps $\pi_\forg$ and $\pi_\quot$ in the $\Gamma_0(N)$ case] \lab{forg-quot-0}
The stack $\sX_0(N)$ (resp.,~$\sX_0(\f{N}{p})$) parametrizes generalized elliptic curves $E \ra S$ equipped with a ``$\Gamma_0(N)$-structure'' (resp.,~a ``$\Gamma_0(\f{N}{p})$-structure''), which on the elliptic curve locus, and for squarefree $N$ also on the entire $S$, is an ample $S$-subgroup $G \subset E^\sm$ of order $N$ (resp.,~$\f{N}{p}$) that is cyclic in the sense of Drinfeld. For general $N$, part of the data of a $\Gamma_0(N)$-structure is a certain open cover $\{S_{(m)}\}_{m \mid N}$ of $S$ over which $G$ is required to live inside suitable ``universal decontractions'' of $E$, see \cite{Ces15a}*{\S5.9} (however, since $p^2 \nmid N$, the $p$-primary part $G[p]$ always lives inside $E$ itself, see \cite{Ces15a}*{5.9.4}). On the elliptic curve locus, the map 
\[
\tst \pi_\forg\colon \sX_0(N) \ra \sX_0(\f{N}{p}) \qq  \text{(resp.,~}\pi_\quot \colon \sX_0(N) \ra \sX_0(\f{N}{p})\text{)}
\]
is defined by replacing $G$ by $G[\f{N}{p}]$ (resp.,~by replacing $E$ by $E/G[p]$ endowed with $G/G[p]$). Granted that for $\pi_\forg$ one contracts $E$ with respect to $G[\f{N}{p}]$, for squarefree $N$ the same definition also works over the entire $\sX_0(N)$, and our task is to explain how to naturally extend it to the entire $\sX_0(N)$ for general $N$. For this, the following lemma suffices because the open substacks $\sX_0(N)_{(m)} \subset \sX_0(N)$ cut out by the $S_{(m)}$ for $m\mid N$ cover $\sX_0(N)$ and pairwise intersect in the elliptic curve locus, to the effect that we only need to define each $(\pi_\forg)|_{\sX_0(N)_{(m)}}$ (resp.,~$(\pi_\quot)|_{\sX_0(N)_{(m)}}$) compatibly with the already given definition on the elliptic curve locus.

For brevity, in the lemma if $E \ra S$ is a generalized elliptic curve with $d$-gon degenerate geometric fibers for some $d \in \bZ_{\ge 1}$, then the \emph{prime to $p$ contraction} of $E$, denoted by $E'$, is the contraction of $E$ with respect to $E^\sm[\f{d}{d^{\ord_p(d)}}]$ (or with respect to any other finite locally free $S$-subgroup that meets the same irreducible components of the geometric fibers of $E \ra S$ as $E^\sm[\f{d}{d^{\ord_p(d)}}]$, see \cite{Ces15a}*{3.2.1}).
\epp

\blem \lab{monster}
In the setting of a prime $p$ and an $N \in \bZ_{\ge 1}$ with $\ord_p(N) = 1$, fix an $m \mid N$, let $E \ra S$ be a generalized elliptic curve equipped with a $\Gamma_0(N)$-structure for which $S_{(m)} = S$, and let $G_{(m)} \subset E^\sm$ be the resulting ample cyclic $S$-subgroup of order $\f{N}{\gcd(m, \f{N}{m})}$ \up{see \cite{Ces15a}*{5.9.4}}.
\benum
\item \lab{mon-a}
There is a unique $\Gamma_0(\f{N}{p})$-structure on $E'$ such that for every fppf $S$-scheme $\wt{S}$ endowed with a generalized elliptic curve $\wt{E} \ra \wt{S}$ that has $m$-gon degenerate geometric fibers and is equipped with an isomorphism between its contraction and $E_{\wt{S}}$, the ample cyclic $S$-subgroups of $\wt{E}'$ of order $\f{N}{p}$ determined by the $\Gamma_0(\f{N}{p})$-structure on $E'$ and by the $\Gamma_0(N)$-structure on $E$ agree.

\item \lab{mon-b}
There is a unique $\Gamma_0(\f{N}{p})$-structure on $E/(G_{(m)}[p])$ such that for every fppf $S$-scheme $\wt{S}$ endowed with an $\wt{E} \ra \wt{S}$ as in \ref{mon-a}, the ample cyclic $S$-subgroup $\wt{G} \subset \wt{E}^\sm$ of order $N$ determined by the $\Gamma_0(N)$-structure on $E$ is such that $\wt{G}/\wt{G}[p] \subset (\wt{E}/\wt{G}[p])^\sm$ agrees with the $S$-subgroup determined by the $\Gamma_0(\f{N}{p})$-structure on $E/(G_{(m)}[p])$.
\eenum
\elem
  
 \bpf 
 An fppf $\wt{S} \ra S$ endowed with a required $\wt{E} \ra \wt{S}$ always exists, see \cite{Ces15a}*{3.2.6}.
 \benum
\item
The uniqueness aspect allows us to work fppf locally on $S$, so we assume that $S = \wt{S}$ and let $\wt{G} \subset \wt{E}^\sm$ be the $S$-subgroup determined by the $\Gamma_0(N)$-structure on $E$. As in \cite{Ces15a}*{\S5.11}, the ample $S$-subgroup $\wt{G}[\f{N}{p}] \subset (\wt{E}')^\sm$ of order $\f{N}{p}$ uniquely extends to a $\Gamma_0(\f{N}{p})$-structure on $E'$. It remains to note that this unique $\Gamma_0(\f{N}{p})$-structure satisfies the claimed compatibility with respect to any other $\wt{E} \ra \wt{S}$ due to \cite{Ces15a}*{5.7}.

\item
To make sense of the characterizing property, one notes that $(G[p]_{(m)})_{\wt{S}}$ identifies with $\wt{G}[p]$ inside $E_{\wt{S}}$ and that $(E/G_{(m)}[p])_{\wt{S}}$ identifies with a contraction of $\wt{E}/\wt{G}[p]$, as is ensured by the uniqueness aspect of \cite{DR73}*{IV.1.2} (see also the review in \cite{Ces15a}*{3.2.1}). 
Granted this, the proof is the same as that of \ref{mon-a} with the role of $\wt{G}[\f{N}{p}]$ replaced by $\wt{G}/\wt{G}[p]$. \qedhere
\eenum
\epf

\bpp[The maps $\pi_\forg$ and $\pi_\quot$ on coarse spaces] \lab{forg-quot-c}
Returning to the setting of \S\ref{setup-2}, we let
\be \lab{fq-coarse}
\pi_\forg,\, \pi_\quot\colon X \rightrightarrows X'
\ee
denote the maps induced on the coarse moduli schemes by the maps $\pi_\forg,\, \pi_\quot\colon \sX \rightrightarrows \sX'$ constructed in \S\S\ref{forg-quot-1}--\ref{forg-quot-0} (we take $N = n$ and $p = 2$). The base changes of the maps \eqref{fq-coarse} to $\bC$ identify with degeneracy maps that appear in a discussion of the theory of newforms, so pullbacks of K\"{a}hler differentials along $(\pi_\forg)_\bQ$ or $(\pi_\quot)_\bQ$ are (associated to) oldforms. The map $\pi_\forg$ induces
\[
\text{an isomorphism} \q X_{\bF_2}^{\mu} \cong X'_{\bF_2} \q \text{and a purely inseparable degree $2$ morphism} \q X_{\bF_2}^{\et} \ra X'_{\bF_2},
\]
as is seen on the elliptic curve locus using \cite{DR73}*{diagram on p.~287}. Analogously, $\pi_\quot$ induces
\[
\text{a purely inseparable degree $2$ morphism} \q X_{\bF_2}^{\mu} \ra X'_{\bF_2} \q \text{and an isomorphism} \q X_{\bF_2}^{\et} \cong X'_{\bF_2}.
\]
\epp

With the maps $\pi_\forg$ and $\pi_\quot$ at our disposal, we turn to producing oldforms that sweep out the $2$-torsion subgroup of $H^0(U^\mu, \Omega^1)/H^0(X_{\bZ_{(2)}}, \Omega)$. This is accomplished in \Cref{lift-2-tor}, with the key step carried out in the following lemma.

\blemt \lab{lift-once}
Let $S \subset X_{\bF_2}^\et$ be the reduced divisor of supersingular points, and let $\Omega|_{X_{\bF_2}^\et}$ be the pullback of the relative dualizing sheaf of $X_{\bZ_{(2)}}$. Every $g \in H^0(X_{\bF_2}^\et, \Omega|_{X_{\bF_2}^\et}(-S))$ lifts to an oldform $\wt{g} \in H^0(X_{\bZ_{(2)}}, \Omega)$ whose restriction to $X_{\bF_2}^\mu$ vanishes.
\elemt

\bpf
Since $X_{\bF_2}$ is semistable (see \S\ref{spec-fib}), we have the $\sO_{X_{\bF_2}^\et}$-module module identification
\be \lab{id-schemes}
\qq \Omega|_{X_{\bF_2}^\et} = \Omega^1_{X_{\bF_2}^\et/\bF_2}(S) \qq \text{inside the generic stalk of $\Omega^1_{X^\et_{\bF_2}/\bF_2}$},
\ee
as may be checked over $\ov{\bF}_2$ by using the theory of regular differentials exposed in \cite{Con00}*{\S5.2}. Therefore, $g$ lies in $H^0(X_{\bF_2}^\et, \Omega^1)$, and hence, via the isomorphism $X_{\bF_2}^\et \cong X'_{\bF_2}$ induced by $\pi_\quot$ (see \S\ref{forg-quot-c}), corresponds to a unique $g' \in H^0(X'_{\bF_2}, \Omega^1)$. By Grothendieck--Serre duality and cohomology and base change (see \cite{Con00}*{Thm.~5.1.2}), we may lift $g'$ to a $G' \in H^0(X'_{\bZ_{(2)}}, \Omega^1)$. We set
\[
\wt{g} \ce (\pi_\quot|_{X_{\bZ_{(2)}} - S})^*(G'|_{X'_{\bZ_{(2)}} - \pi_\quot(S)}) \in H^0(X_{\bZ_{(2)}} - S, \Omega^1) = H^0(X_{\bZ_{(2)}}, \Omega)
\]
(the normality of $X_{\bZ_{(2)}}$ ensures the equality because $\Omega$ is a line bundle whose restriction to $X_{\bZ_{(2)}} - S$ is $\Omega^1$). By construction, $\wt{g}$ is an oldform in $H^0(X_{\bZ_{(2)}}, \Omega)$ that agrees with $g$ on $X_{\bF_2}^\et$ and that vanishes on $X_{\bF_2}^\mu$ because the map $X_{\bF_2}^\mu \ra X'_{\bF_2}$ induced by $\pi_\quot$ is purely inseparable of degree $2$.
\epf

\bpropt \lab{lift-2-tor}
Every element of $\p{H^0(U^\mu, \Omega^1)/H^0(X_{\bZ_{(2)}}, \Omega)}[2]$ lifts to an oldform in $H^0(U^\mu, \Omega^1)$.
\epropt

\bpf
The stalks of $\sO_{X_{\bZ_{(2)}}}$ at the generic points of $X_{\bF_2}^\mu$ and $X_{\bF_2}^\et$ are discrete valuation rings with $2$ as a uniformizer.  Thus, as explained in \cite{BLR90}*{p.~104}, there are ``order functions'' $v^\mu$ and $v^\et$ that measure the valuations of sections of $\Omega$ at these respective points. In this notation, $H^0(U^\mu, \Omega^1)$ (resp.,~$H^0(X_{\bZ_{(2)}}, \Omega)$) identifies with the set of $f \in H^0(X_\bQ, \Omega^1)$ for which $v^\mu(f) \ge 0$ (resp.,~for which $v^\mu(f) \ge 0$ and $v^\et(f) \ge 0$), similarly to the proof of \Cref{hop-f}.

If $f \in H^0(U^\mu, \Omega^1)$ represents an element of $(H^0(U^\mu, \Omega^1)/H^0(X_{\bZ_{(2)}}, \Omega))[2]$, then $2f$ is a global section of $\Omega$ that vanishes on $X_{\bF_2}^\mu$. The restriction of $2f$ to $X_{\bF_2}^\et$ therefore lies in $H^0(X_{\bF_2}^\et, \Omega|_{X_{\bF_2}^\et}(-S))$, as is required for \Cref{lift-once} to apply. \Cref{lift-once} supplies an oldform $\wt{g} \in H^0(X_{\bZ_{(2)}}, \Omega)$ that agrees with $2f$ on $X_{\bF_2}^\et$ and vanishes on $X_{\bF_2}^\mu$. It remains to note that, by the discussion of the previous paragraph, $\f{\wt{g}}{2}$ is an oldform in $H^0(U^\mu, \Omega^1)$ for which 
\[
\tst \wt{g}- 2f \in 2\cdot H^0(X_{\bZ_{(2)}}, \Omega), \qq \text{i.e.,} \qq \f{\wt{g}}{2} - f \in H^0(X_{\bZ_{(2)}}, \Omega). \qedhere
\]
\epf

\bremt \lab{odd-reproof}
In the case $\Gamma = \Gamma_0(n)$, the discussion of \S\ref{forg-quot-c} and the proofs of \Cref{lift-once} and \Cref{lift-2-tor} are not specific to the prime $p = 2$ (see also Remark \ref{also-odd-p}). In particular, if $n \in \bZ_{\ge 1}$ and $p$ is a prime with $\ord_p(n) \le 1$, then, adopting the analogous notation $U^\mu \subset X_0(n)_{\bZ_{(p)}}$ also for odd $p$, they show that every $p$-torsion element of $H^0(U^\mu, \Omega^1)/H^0(X_0(n)_{\bZ_{(p)}}, \Omega)$ lifts to an oldform in $H^0(U^\mu, \Omega^1)$. Since for odd $p$ every element of this quotient is $p$-torsion by \cite{Edi06}*{2.5 and 2.7} (equivalently, by the proof of \Cref{int-thy} below), this combines with \Cref{main-reduction} and \Cref{more-odd-p} to reprove that the Manin constant of a new elliptic optimal quotient of $J_0(n)$ is not divisible by any odd prime $p$ with $\ord_p(n) \le 1$. The key distinction of this proof is that it does not use exactness results for semiabelian N\'{e}ron models over bases of low ramification (compare with the well-known argument recalled in the proof of \Cref{odd-p}).
\eremt

For extending \Cref{lift-2-tor} beyond the $2$-torsion subgroup of $H^0(U^\mu, \Omega^1)/H^0(X_{\bZ_{(2)}}, \Omega)$, it will be handy to work with the Deligne--Mumford stack $\sX_{\bZ_{(2)}}$ instead of its coarse moduli scheme $X_{\bZ_{(2)}}$. To facilitate for this, we supplement the discussion of \S\ref{spec-fib} with a similar discussion of $\sX_{\bZ_{(2)}}$. The principal advantage of $\sX_{\bZ_{(2)}}$ over $X_{\bZ_{(2)}}$ is its regularity (see \S\ref{setup-2}), which will permit an effective use of techniques from intersection theory. The principal disadvantage is the loss of Grothendieck--Serre duality, which was an important component of the proof of \Cref{lift-once}.

\bpp[The structure of $\sX_{\bZ_{(2)}}$] \lab{stack-str}
Due to the moduli interpretation of $\sX_{\bZ_{(2)}}$, the map 
\[
\sX_{\bZ_{(2)}} \ra X_{\bZ_{(2)}}
\]
towards the coarse moduli scheme is \'{e}tale over the locus of $X_{\bZ_{(2)}}$ on which the $j$-invariant satisfies $j \neq 0$ and $j \neq 1728$, see \cite{Ces15a}*{4.1.4 and proof of Thm.~6.7}. Thus, since $\sX$ is normal, $\sX_{\bZ_{(2)}} \ra \Spec \bZ_{(2)}$ is smooth away from the supersingular points of its $\bF_2$-fiber $\sX_{\bF_2}$ (see \S\ref{spec-fib}). In particular, by the (R$_0$)$+$(S$_1$) criterion, $\sX_{\bF_2}$ is reduced and, by \cite{DR73}*{VI.6.10}, $X_{\bF_2}$ is the coarse moduli space of $\sX_{\bF_2}$. 

In contrast, the smooth locus of $\sX'_{\bZ_{(2)}} \ra \Spec \bZ_{(2)}$ is the entire $\sX'_{\bZ_{(2)}}$ by \cite{DR73}*{IV.6.7}.

The decomposition of $X_{\bF_2}$ into irreducible components $X_{\bF_2}^\mu$ and $X_{\bF_2}^\et$ corresponds to the decomposition of $\sX_{\bF_2}$ into irreducible components 
\[
\sX_{\bF_2}^\mu \q \text{and} \q \sX_{\bF_2}^\et.
\]
We let $\sS \subset \sX_{\bF_2}$ denote the reduced closed substack consisting of the supersingular points, so that $\sX_{\bF_2}^\mu$ and $\sX_{\bF_2}^\et$ meet precisely at $\sS$. Due to \cite{DR73}*{V.1.16 (ii)}, the intersections of $\sX_{\bF_2}^\mu$ and $\sX_{\bF_2}^\et$ in $\sX_{\bF_2}$ at the points of $\sS$ are transversal. Thus, due to the regularity of $\sX_{\bZ_{(2)}}$, the intersections of $\sX_{\bF_2}^\mu$ and $\sX_{\bF_2}^\et$ in $\sX_{\bZ_{(2)}}$ are also transversal. 

By \cite{DR73}*{V.1.16 (i)}, if we let the replacement of $\sX$ by $\sY$ indicate the elliptic curve locus, then we have the isomorphism
\be \lab{mu-et-pr}
\sY'_{\bF_2} \isomto \sY^\mu_{\bF_2} \qq \text{(resp.,~}\sY'_{\bF_2} \isomto \sY^\et_{\bF_2}) 
\ee
obtained by supplementing the universal elliptic curve of $\sY'_{\bF_2}$ (resp.,~the Frobenius pullback of the universal elliptic curve of $\sY'_{\bF_2}$) with the subgroup of order $2$ given by the kernel of Frobenius (resp.,~the kernel of Verschiebung). 
\epp

\bpp[The line bundle $\omega$] \lab{omega}
The cotangent space at the identity section of the universal generalized elliptic curve gives a line bundle $\omega$ on $\sX$ (resp.,~on $\sX'$) of weight $1$ modular forms (we sometimes write $\omega_\sX$, etc.~to emphasize the space on which $\omega$ lives). We will primarily be concerned with cusp forms, so we let `$\cusps$' denote the reduced relative effective Cartier divisor on $\sX$ (resp.,~on $\sX'$) over $\bZ$ cut out by the degeneracy locus of the universal generalized elliptic curve (see \cite{Ces15a}*{4.4.2~(b) and 5.13~(b)}); \emph{a posteriori}, `$\cusps$' is also the reduced complement of the elliptic curve locus of $\sX$ or $\sX'$. Depending on the context, we also write `$\cusps$' for base changes or restrictions of this divisor. The line bundle whose global sections are weight $2$ cusp forms is therefore $\omega^{\tensor 2}(-\cusps)$. 

The modular definitions of the maps $\pi_\forg$ and $\pi_\quot$ given in \S\S\ref{forg-quot-1}--\ref{forg-quot-0} also produce underlying $\sX$-morphisms from the universal generalized elliptic curve of $\sX$ to the pullback of the universal generalized elliptic curve of $\sX'$. The effect of these $\sX$-morphisms on the cotangent spaces at the identity sections gives rise to $\sO_\sX$-module morphisms
\be \lab{omega-pull}
\pi_\forg^*(\omega_{\sX'}) \ra \omega_{\sX} \qq \text{and} \qq \pi_\quot^*(\omega_{\sX'}) \ra \omega_{\sX}.
\ee
Thus, since $\pi_\forg$ and $\pi_\quot$ are finite locally free (see \S\ref{forg-quot-1}) and when restricted to `$\cusps$' on $\sX$ factor through `$\cusps$' on $\sX'$, the morphisms \eqref{omega-pull} lead to $\sO_{\sX}$-module morphisms
\[
\pi_\forg^*(\omega^{\tensor 2}_{\sX'}(-\cusps)) \ra \omega^{\tensor 2}_{\sX}(-\cusps) \qq \text{and} \qq \pi_\quot^*(\omega^{\tensor 2}_{\sX'}(-\cusps)) \ra \omega^{\tensor 2}_{\sX}(-\cusps).
\]
We analyze the restrictions of these morphisms to $\sX_{\bF_2}^\mu$ and $\sX_{\bF_2}^\et$ in the following lemma.
\epp

\blemt \lab{pi-fest} \hfill
\benum
\item \lab{PF-a}
The map $\pi_\forg$ restricts to an isomorphism $\pi_\forg\colon \sX_{\bF_2}^\mu \isomto \sX'_{\bF_2}$ for which the pullback
\[
\pi_\forg^*(\omega_{\sX'_{\bF_2}}^{\tensor 2}(-\cusps)) \ra \omega^{\tensor 2}_{\sX_{\bF_2}^\mu}(-\cusps)
\]
is also an isomorphism.

\item \lab{PF-c}
The map $\pi_\quot$ restricts to an isomorphism $\pi_\quot\colon \sX_{\bF_2}^\et \isomto \sX'_{\bF_2}$ for which the pullback
\[
\pi_\quot^*(\omega_{\sX'_{\bF_2}}^{\tensor 2}(-\cusps)) \ra \omega^{\tensor 2}_{\sX_{\bF_2}^\et}(-\cusps)
\]
induces an identification $\pi_\quot^*(\omega_{\sX'_{\bF_2}}^{\tensor 2}(-\cusps)) \cong \omega^{\tensor 2}_{\sX_{\bF_2}^\et}(-\cusps - 2\sS)$.

\item \lab{PF-b}
The map $\pi_\quot$ restricts to a morphism $\pi_\quot\colon \sX_{\bF_2}^\mu \ra \sX'_{\bF_2}$ for which the pullback
\[
\pi_\quot^*(\omega_{\sX'_{\bF_2}}^{\tensor 2}(-\cusps)) \ra \omega^{\tensor 2}_{\sX_{\bF_2}^\mu}(-\cusps)
\]
vanishes.
\eenum
\elemt

\bpf 
The divisor `$\cusps$' on $\sX'_{\bZ_{(2)}}$ is \'{e}tale over $\bZ_{(2)}$, and hence is also $\bZ_{(2)}$-fiberwise reduced---this follows from \cite{Ces15a}*{4.4.2 (b) and 5.13 (b)}.
\benum
\item
On the elliptic curve locus the claim follows from the description of the first map of \eqref{mu-et-pr}. Thus, $\pi_\forg\colon \sX_{\bF_2}^\mu \ra \sX'_{\bF_2}$ is finite locally free of rank $1$ (see \S\ref{forg-quot-1}), and hence is an isomorphism. Moreover, since its restriction to `$\cusps$' of the source factors through `$\cusps$' of the target, the reducedness of the latter ensures that $\pi_\forg\colon \sX_{\bF_2}^\mu \isomto \sX'_{\bF_2}$ identifies `$\cusps$' of its source and target. It remains to note that $\pi_\forg^*(\omega_{\sX'_{\bF_2}}^{\tensor 2}) \isomto \omega^{\tensor 2}_{\sX_{\bF_2}^\mu}$ because $\pi_\forg$ does not change the relative identity component of the smooth locus of the universal generalized elliptic curve.

\item
As in the proof of \ref{PF-a}, the map $\pi_\quot \colon \sX_{\bF_2}^\et \ra \sX'_{\bF_2}$ is an isomorphism that identifies `$\cusps$' of its source and target. It remains to show that the pullback map
\be \lab{Hasse-pull}
\q \pi_\quot^*(\omega_{\sX'_{\bF_2}}) \ra \omega_{\sX_{\bF_2}^\et}
\ee
induces an identification 
\[
\qq \pi_\quot^*(\omega_{\sX'_{\bF_2}}) \cong \omega_{\sX_{\bF_2}^\et}(-\sS).
\]
The subgroup of order $2$ coming from the $2$-primary part of the $\Gamma$-structure on the universal generalized elliptic curve of $\sX_{\bF_2}^\et - \sS$ is \'{e}tale because the locus where this subgroup is of multiplicative type is \emph{a priori} open and does not meet the elliptic curve locus (and hence is empty). Therefore, \eqref{Hasse-pull}, being induced by pullback along the quotient by this $2$-primary part, is an isomorphism away from the supersingular points.

For the remaining claim that the divisor cut out by \eqref{Hasse-pull} is precisely $\sS$, we may work on the elliptic curve locus and after restriction along the isomorphism $\sY'_{\bF_2} \isomto \sY_{\bF_2}^\et$ of \eqref{mu-et-pr}. After this restriction, \eqref{Hasse-pull} identifies with the map induced by pullback along the Verschiebung isogeny of the universal elliptic curve of $\sY'_{\bF_2}$, i.e.,~with the Hasse invariant. It remains to recall from \cite{KM85}*{12.4.4} that the Hasse invariant has simple zeroes at the supersingular~points.

\item
For the sought vanishing we may work on the elliptic curve locus, so, due to \eqref{mu-et-pr}, it suffices to note that for an elliptic curve over a base scheme of characteristic $p > 0$ the pullback of any differential form along the relative Frobenius morphism vanishes.
\qedhere
\eenum
\epf

We are ready for the following key proposition, which will replace \Cref{lift-once} when attempting to lift arbitrary elements of $H^0(U^\mu, \Omega^1)/H^0(X_{\bZ_{(2)}}, \Omega)$ to oldforms. The role of its surjectivity assumption is to serve as a replacement of the Grothendieck--Serre duality for the Deligne--Mumford stack $\sX'_{\bZ_{(2)}}$ (such duality for the scheme $X'_{\bZ_{(2)}}$ was important in the proof of \Cref{lift-once}; see also \Cref{KS}~\ref{KS-a}).

\bpropt \lab{lift-twice}
If the pullback map 
\[
H^0(\sX'_{\bZ_{(2)}}, \omega^{\tensor 2}(-\cusps)) \ra H^0(\sX'_{\bF_2}, \omega^{\tensor 2}(-\cusps))
\]
is surjective, then every $g \in H^0(\sX_{\bF_2}^\mu, \omega^{\tensor 2}(-\cusps - \sS))$ lifts to an oldform $\wt{g} \in H^0(\sX_{\bZ_{(2)}}, \omega^{\tensor 2}(-\cusps))$ that vanishes on $\sX_{\bF_2}^\et$.
\epropt

\bpf
By \Cref{pi-fest}~\ref{PF-a}, $g$ is the $\pi_\forg$-pullback of a unique $g' \in H^0(\sX'_{\bF_2}, \omega^{\tensor 2}(-\cusps - \sS'))$, where $\sS' \subset \sX'_{\bF_{2}}$ is the reduced closed substack supported at the supersingular points. Due to the surjectivity assumption, $g'$ lifts to a $G' \in H^0(\sX'_{\bZ_{(2)}}, \omega^{\tensor 2}(-\cusps))$. We set
\[
\wt{g}_0 \ce \pi_\forg^*(G') \in H^0(\sX_{\bZ_{(2)}}, \omega^{\tensor 2}(-\cusps)),
\]
so that $\wt{g}_0$ is an oldform that lifts $g$. We claim that the restriction 
\be \lab{lies-in}
h \ce \wt{g}_0|_{\sX_{\bF_2}^\et} \in H^0(\sX_{\bF_2}^\et, \omega^{\tensor 2}(-\cusps)) \qq \text{lies in} \qq H^0(\sX_{\bF_2}^\et, \omega^{\tensor 2}(-\cusps - 2\sS)).
\ee
For this, we may work on the elliptic curve locus and after restricting along the isomorphism $\sY'_{\bF_2} \isomto \sY_{\bF_2}^\et$ of \eqref{mu-et-pr}. Under this isomorphism, $\omega^{\tensor 2}_{\sY_{\bF_2}^\et}$ identifies with $\omega^{\tensor 4}_{\sY'_{\bF_2}}$ and $\pi_\forg|_{\sY_{\bF_2}^\et}$ identifies with the Frobenius morphism of $\sY'_{\bF_2}$, so $h$ identifies with the Frobenius pullback of $g'$. Since $g'$, when viewed as a global section of $\omega^{\tensor 2}_{\sX'_{\bF_2}}$, vanishes on $\sS'$, \eqref{lies-in} follows from the fact that the Frobenius pullback of $\sO_{\sY'_{\bF_2}}(-\sS')$ is $\sO_{\sY'_{\bF_2}}(-2\sS')$.

Due to \eqref{lies-in} and \Cref{pi-fest}~\ref{PF-c}, $h$ is the $\pi_\quot$-pullback of a unique $h' \in H^0(\sX'_{\bF_2}, \omega^{\tensor 2}(-\cusps))$. The surjectivity assumption lifts $h'$ to an $H' \in H^0(\sX'_{\bZ_{(2)}}, \omega^{\tensor 2}(-\cusps))$ and we set
\[
\wt{h} \ce \pi_\quot^*(H') \in H^0(\sX_{\bZ_{(2)}}, \omega^{\tensor 2}(-\cusps)).
\]
By construction and \Cref{pi-fest}~\ref{PF-b}, $\wt{h}$ is an oldform, agrees with $h$ on $\sX^\et_{\bF_2}$, and vanishes on $\sX^\mu_{\bF_2}$. In conclusion, the oldform $\wt{g} \ce \wt{g}_0 - \wt{h} \in H^0(\sX_{\bZ_{(2)}}, \omega^{\tensor 2}(-\cusps))$ lifts $g$ and vanishes on $\sX_{\bF_2}^\et$.
\epf

In order to take advantage of \Cref{lift-twice}, we seek to relate $\omega^{\tensor 2}(-\cusps)$ to the ``relative dualizing sheaf'' $\Omega$ of \S\ref{rel-dual} via Kodaira--Spencer type isomorphisms supplied by the following lemma.

\blemt \lab{KS} \hfill
\benum
\item \lab{KS-a}
On $\sX'_{\bZ_{(2)}}$, there is an $\sO_{\sX'_{\bZ_{(2)}}}$-module isomorphism $\Omega^1 \cong \omega^{\tensor 2}(-\cusps)$.

\item \lab{KS-b}
On $\sX_{\bZ_{(2)}}$, there is an $\sO_{\sX_{\bZ_{(2)}}}$-module isomorphism $\Omega \cong \omega^{\tensor 2}(-\cusps+ \sX_{\bF_2}^\et)$.
\eenum
\elemt

\bpf
We will bootstrap the claims from their analogue for $\sX(1)$ supplied by \cite{Kat73}*{A1.3.17}:
\be \lab{Katz-input}
\Omega^1_{\sX(1)/\bZ} \cong \omega^{\tensor 2}_{\sX(1)}(-\cusps).
\ee 
For any congruence subgroup $H \subset \GL_2(\wh{\bZ})$, the structure map $\pi\colon \sX_H \ra \sX(1)$ is finite locally free, so, due to the base change compatibility of the formation of the relative dualizing sheaf (see \cite{Con00}*{Thm.~3.6.1}), there is a ``relative dualizing'' $\sO_{\sX_H}$-module $\Omega_{\sX_H/\sX(1)}$ constructed \'{e}tale locally on $\sX(1)$. Explicitly, due to \cite{Con00}*{bottom half of p.~31 and pp.~137--139, esp.~(VAR6) on p.~139, supplemented by Cor.~3.6.4}, $\Omega_{\sX_H/\sX(1)}$ identifies with $\sH om_{\sO_{\sX(1)}}(\pi_*(\sO_{\sX_H}), \sO_{\sX(1)})$ regarded as an $\sO_{\sX_H}$-module.

By working \'{e}tale locally on $\sX(1)$, \cite{Con00}*{Thm.~4.3.3 and (4.3.7); see also bottom of p.~206} supply an $\sO_{\sX_H}$-module isomorphism
\be \lab{duality-at-its-best}
\Omega_{\sX_H/\sX(1)} \otimes_{\sO_{\sX_H}} \pi^* \Omega_{\sX(1)/\bZ} \cong \Omega_{\sX_H/\bZ}.
\ee
To proceed further, we assume that $\sX_H$ is regular, so that $\pi$ is a local complete intersection (see \cite{Liu02}*{6.3.18}), and hence has Gorenstein fibers, to the effect that $\Omega_{\sX_H/\sX(1)}$ is a line bundle (see \cite{Con00}*{Thm.~3.5.1}). Then, since $\pi$ is \'{e}tale over a dense open of $\sX(1)$, the element 
\[
\mathrm{trace} \in \Hom_{\sO_{\sX(1)}}(\pi_*(\sO_{\sX_H}), \sO_{\sX(1)}) \cong \Gamma(\sX_H, \Omega_{\sX_H/\sX(1)})
\]
gives rise to the identification
\be \lab{omega-kick}
\tst \Omega_{\sX_H/\sX(1)}\cong \sO_{\sX_H}(\sum_{x \in \abs{\sX_H}^{(1)}} d_x \cdot \ov{\{x\}}),
\ee
where the sum runs over the height $1$ points $x$ of $\sX_H$, the corresponding to $x$ irreducible Weil divisor on $\sX_H$ is denoted by $\ov{\{ x\}}$, and $d_x$ denotes the valuation of the different ideal of the extension $\sO_{\sX_H, x}^\sh/\sO_{\sX(1), \pi(x)}^\sh$ of discrete valuation rings. Since $d_x = 0$ whenever this extension is \'{e}tale, each $x$ that contributes to the sum either is the generic point of an irreducible component of a closed fiber of $\sX_H \ra \Spec \bZ$ or lies on the cusps of $(\sX_H)_\bQ$. Moreover, at the latter $x$ the ramification is tame and $d_x = e_x - 1$, where $e_x$ is the ramification index of $\sO_{\sX_H, x}^\sh/\sO_{\sX(1), \pi(x)}^\sh$. By combining this with \eqref{Katz-input}--\eqref{omega-kick}, we arrive at the identification
\be \lab{omg-omg}
\tst \Omega_{\sX_H/\bZ} \cong (\pi^*(\omega_{\sX(1)}))^{\otimes 2}(-\cusps + \sum_y d_y \cdot \ov{\{ y \}}),
\ee
where $y$ runs over the generic points of the irreducible components of the closed $\bZ$-fibers of $\sX_H$ and `cusps' denotes the reduced complement of the elliptic curve locus of $\sX_H$. 

In the case when $\sX_H$ is $\sX$ or $\sX'$, the map $\pi$ is the forgetful contraction and does not change the relative identity component of the smooth locus of the universal generalized elliptic curve, so $\pi^*(\omega_{\sX(1)})$ identifies with $\omega_\sX$ or $\omega_{\sX'}$, respectively. Therefore, since $\Omega_{\sX'_{\bZ_{(2)}}/\bZ_{(2)}} \cong \Omega_{\sX'_{\bZ_{(2)}}/\bZ_{(2)}}^1$ due to the ${\bZ_{(2)}}$-smoothness of $\sX'_{\bZ_{(2)}}$, the sought conclusion will follow from \eqref{omg-omg} once we identify the $d_y$ for $y$ of residue characteristic $2$ in the case when $\sX_H = \sX$ or $\sX_H = \sX'$.
\benum
\item
Since the ``level'' $\f{n}{2}$ of $\sX'$ is odd, the map $\sX'_{\bZ_{(2)}} \ra \sX(1)_{\bZ_{(2)}}$ is \'{e}tale over a fiberwise dense open of $\sX(1)_{\bZ_{(2)}}$. Therefore, $d_y = 0$ whenever $y$ has residue characteristic $2$.

\item
By \Cref{pi-fest}~\ref{PF-a} and the proof of \ref{KS-a}, $\pi$ is generically \'{e}tale on $\sX_{\bF_2}^\mu$. In contrast, \eqref{mu-et-pr} identifies the map $\sY_{\bF_2}^\et \ra \sY'_{\bF_2}$ induced by $\pi_\forg$ with the Frobenius of $\sY'_{\bF_2}$, which is not generically \'{e}tale. We therefore conclude from \eqref{omg-omg} and \cite{AS02}*{A.3} that
\be \lab{omg-almost}
\qq \Omega_{\sX_{\bZ_{(2)}}/\bZ_{(2)}} \cong \omega^{\tensor 2}(-\cusps + d\cdot \sX_{\bF_2}^\et) \qq \text{for some} \q  d \ge 1.
\ee
To see that also $d \le 1$, we note that on strict Henselizations of $\sX$ and $\sX'$ at the generic points of $\sX_{\bF_2}^\et$ and $\sX_{\bF_2}'$ the map $\pi_\forg$ induces a degree $2$ extension of absolutely unramified discrete valuation rings (in particular, the trace of $1$ in this extension is a uniformizer).
\qedhere
\eenum
\epf

\bremt \lab{KS-old}
The isomorphism $H^0(\sX_\bQ, \Omega^1) \cong H^0(\sX_\bQ, \omega^{\tensor 2}(-\cusps))$ supplied by \Cref{KS}~\ref{KS-b} or even by any $\sO_{\sX_\bQ}$-module isomorphism $\Omega^1 \simeq \omega^{\tensor 2}(-\cusps)$, any two of which differ by $\bQ^\times$-scaling, identifies the spaces of oldforms on both sides, as may be checked over $\bC$ (see also \cite{Gro90}*{3.15}).
\eremt

To proceed beyond the $2$-torsion subgroup treated in \Cref{lift-2-tor}, we begin by recording the following basic fact about the structure of $H^0(U^\mu, \Omega^1)/H^0(X_{\bZ_{(2)}}, \Omega)$. Its proof below is modeled on that of \cite{Edi06}*{2.7} given there in the $\Gamma_0(n)$ context (see also \cite{DR73}*{VII.3.19--20} and \cite{BDP16}*{B.3.2} for similar results). The use of the Atkin--Lehner involution in the proof is for convenience and compensates for the fact that \Cref{pi-fest}~\ref{PF-a} is specific to the irreducible component $\sX_{\bF_2}^\mu$ of $\sX_{\bF_2}$.

\bpropt \lab{int-thy}
The quotient $H^0(U^\mu, \Omega^1)/H^0(X_{\bZ_{(2)}}, \Omega)$ is killed by $4$.
\epropt

\bpf
As noted towards the end of \S\ref{spec-fib}, both $H^0(U^\mu, \Omega^1)$ and $H^0(X_{\bZ_{(2)}}, \Omega)$ are $\bZ_{(2)}$-lattices inside $H^0(X_\bQ, \Omega^1)$. Therefore, every element of their quotient is killed by a power of $2$. 

Due to the moduli interpretation of $\sX$ (combined with \cite{KM85}*{6.1.1 (1)} in the $\Gamma_1(n)$-case), the Atkin--Lehner involution $w_2$ of $X_\bQ$ extends (uniquely) to an involution of the elliptic curve locus of $X_{\bZ_{(2)}}$, and this extension interchanges the generic points of the irreducible components of $X_{\bF_2}$. Therefore, the automorphism of $H^0(X_\bQ, \Omega^1)$ induced by $w_2$ respects the $\bZ_{(2)}$-lattice $H^0(X_{\bZ_{(2)}}, \Omega)$ and interchanges $H^0(U^\mu, \Omega^1)$ and $H^0(U^\et, \Omega^1)$ (see the first paragraph of the proof of \Cref{lift-2-tor}). Our task becomes showing that $H^0(U^\et, \Omega^1)/H^0(X_{\bZ_{(2)}}, \Omega)$ is killed by $4$.

By \Cref{app-main}~\ref{AM-b} and Remark \ref{AM-X0n} (see also \S\ref{stack-str}), we have compatible identifications
\be \lab{U-U-id}
\qq H^0(U^\et, \Omega^1) \cong H^0(\sU^\et, \Omega^1) \qq \text{and} \qq H^0(X_{\bZ_{(2)}}, \Omega) \cong H^0(\sX_{\bZ_{(2)}}, \Omega),
\ee
so we may switch to working with stacks. The principal advantage in this is that due to its regularity, $\sX_{\bZ_{(2)}}$ admits a robust intersection theory formalism (see, for instance, \cite{BDP16}*{\S B.2.2}) analogous to the case of a proper regular arithmetic surface.

We fix an $f \in H^0(\sU^\et, \Omega^1) \setminus H^0(\sX_{\bZ_{(2)}}, \Omega)$, let $m > 0$ be minimal such that $2^m f \in H^0(\sX_{\bZ_{(2)}}, \Omega)$, and seek to show that $m \le 2$ by using the fact that $2^mf$ does not vanish on $\sX_{\bF_2}^\mu$ but vanishes to order at least $m$ along $\sX_{\bF_2}^\et$.  Since the intersections of $\sX_{\bF_2}^\mu$ and $\sX_{\bF_2}^\et$ in $\sX_{\bZ_{(2)}}$ are transversal, the restriction of $2^mf$ to $\sX_{\bF_2}^\mu$ identifies with a nonzero global section of the line bundle
\[
\Omega|_{\sX_{\bF_2}^\mu}(-m\sS) \overset{\text{\ref{KS}~\ref{KS-b}}}{\cong} \omega^{\tensor 2}_{\sX_{\bF_2}^\mu}(- \cusps + (1 - m)\sS) \overset{\text{\ref{pi-fest}~\ref{PF-a}}}{\cong} \omega^{\tensor 2}_{\sX_{\bF_2}'}(- \cusps + (1 - m)\sS'),
\]
whose degree must therefore be nonnegative (we let $\sS'$ be the image of $\sS$ under $\pi_\forg\colon \sX_{\bF_2}^\mu \isomto \sX_{\bF_2}'$). The sought $m \le 2$ follows by taking into account the isomorphism $\omega_{\sX_{\bF_2}'} \cong \sO_{\sX_{\bF_2}'}(\sS')$ supplied by the Hasse invariant (see the proof of \Cref{pi-fest}~\ref{PF-c}) and by recalling that `$\cusps$' $\neq \emptyset$.
\epf

\bremt \lab{Raynaud}
By \Cref{lift-2-tor}, every $2$-torsion element of $H^0(U^\mu, \Omega^1)/H^0(X_{\bZ_{(2)}}, \Omega)$ lifts to an oldform in $H^0(U^\mu, \Omega^1)$, so \Cref{int-thy} shows that the finite cyclic group 
\[
\tst \f{H^0(U^\mu, \Omega^1)}{H^0(X_{\bZ_{(2)}}, \Omega) + H^0(U^\mu, \Omega^1) \cap (\bQ \cdot f)^\perp}
\]
that appears in \Cref{main-reduction} is killed by $2$. Therefore, \Cref{odd-n,main-reduction} reprove a result of Mazur--Raynaud, \cite{AU96}*{Prop.~3.1}: in the setting of \S\ref{setup}, if $\ord_2(n) \le 1$, then the Manin--Stevens constant $c_\pi$ satisfies $\ord_2(c_\pi) \le 1$. Similarly to \Cref{odd-reproof}, the distinction of this reproof is that it does not use exactness results for semiabelian N\'{e}ron models.
\eremt

We are ready to investigate the liftability to oldforms in $H^0(U^\mu, \Omega^1)$ of arbitrary elements of  $H^0(U^\mu, \Omega^1)/H^0(X_{\bZ_{(2)}}, \Omega)$ (\Cref{lift-2-tor} only addressed elements killed by $2$).

\bthmt \lab{lift-all}
If the pullback map 
\[
H^0(\sX'_{\bZ_{(2)}}, \Omega^1) \ra H^0(\sX'_{\bF_2}, \Omega^1)
\]
is surjective, then every element of $H^0(U^\mu, \Omega^1)/H^0(X_{\bZ_{(2)}}, \Omega)$ lifts to an oldform in $H^0(U^\mu, \Omega^1)$.
\ethmt

\bpf
Similarly to the proof of \Cref{int-thy}, the Atkin--Lehner involution $w_2$ reduces us to showing that every element of 
\[
H^0(U^\et, \Omega^1)/H^0(X_{\bZ_{(2)}}, \Omega)
\]
lifts to an oldform in $H^0(U^\et, \Omega^1)$, and we already know such liftability for $2$-torsion elements due to \Cref{lift-2-tor}. Moreover, the identifications \eqref{U-U-id} permit us to switch to working with stacks. In conclusion, we seek to show that for every 
\be \lab{f-0}
f_0 \in H^0(\sU^\et, \Omega^1) \qq \text{such that} \qq 2f_0 \not \in H^0(\sX_{\bZ_{(2)}}, \Omega)
\ee
there exists some oldform $\wt{f}_0 \in H^0(\sU^\et, \Omega^1)$ for which $2(f_0 - \wt{f}_0) \in H^0(\sX_{\bZ_{(2)}}, \Omega)$.

We set $f \ce 4f_0$, so that, by \Cref{int-thy} and \eqref{f-0}, $f$ is a global section of $\Omega$ on $\sX_{\bZ_{(2)}}$ that vanishes to order $2$ along $\sX_{\bF_2}^\et$ but does not vanish on $\sX_{\bF_2}^\mu$. In particular, under the isomorphism
\be \lab{KS-input}
\Omega \overset{\text{\ref{KS} \ref{KS-b}}}{\cong} \omega^{\tensor 2}(-\cusps + \sX_{\bF_2}^\et),
\ee
$f$ lies in $\omega^{\tensor 2}(-\cusps)$ and vanishes on $\sX_{\bF_2}^\et$, so its pullback to $H^0(\sX_{\bF_2}^\mu, \omega^{\tensor 2}(-\cusps))$ lies in $H^0(\sX_{\bF_2}^\mu, \omega^{\tensor 2}(-\cusps - \sS))$. Therefore, \Cref{lift-twice} (with \Cref{KS}~\ref{KS-a}) supplies an oldform 
\[
\wt{f} \in H^0(\sX_{\bZ_{(2)}}, \omega^{\tensor 2}(-\cusps))
\]
that agrees with $f$ on $\sX_{\bF_2}^\mu$ and vanishes on $\sX_{\bF_2}^\et$. This $\wt{f}$ satisfies $f - \wt{f} \in 2 \cdot H^0(\sX_{\bZ_{(2)}}, \omega^{\tensor 2}(-\cusps))$ and, when viewed as a global section of $\Omega$ via \eqref{KS-input}, is an oldform (see \Cref{KS-old}) that vanishes to order at least $2$ along $\sX_{\bF_2}^\et$. The oldform $\f{\wt{f}}{4}$ is then a sought $\wt{f}_0$.
\epf

The following lemma helps us recognize situations in which \Cref{lift-all} applies, i.e.,~in which the surjectivity assumption holds.

\blemt \lab{surj-crit} Fix an odd $m \in \bZ_{\ge 1}$.
\benum
\item \lab{SC-a}
The pullback map 
\[
\qq H^0(\sX_1(m)_{\bZ_{(2)}}, \Omega^1) \ra H^0(\sX_1(m)_{\bF_2}, \Omega^1)
\]
is surjective whenever $\sX_1(m)_{\bZ_{(2)}}$ is a scheme, for instance, whenever $m > 3$.

\item \lab{SC-b}
The pullback map 
\[
\qq H^0(\sX_0(m)_{\bZ_{(2)}}, \Omega^1) \ra H^0(\sX_0(m)_{\bF_2}, \Omega^1)
\]
is surjective if and only if so is the pullback map 
\be \lab{coarse-thing}
\qq H^0(X_0(m)_{\bF_2}, \Omega^1) \ra H^0(\sX_0(m)_{\bF_2}, \Omega^1),
\ee
and this is the case if $m$ is a prime or if $m$ has a prime factor $q$ with $q \equiv 3 \bmod 4$.

\eenum
\elemt

\bpf \hfill
\benum
\item
Since $m$ is odd, $\sX_1(m)_{\bZ_{(2)}} \ra \Spec \bZ_{(2)}$ is proper and smooth of relative dimension $1$ (see \cite{DR73}*{IV.6.7}). Therefore, if $\sX_1(m)_{\bZ_{(2)}}$ is a scheme, then the surjectivity in question follows from the formalism of Grothendieck--Serre duality and cohomology and base change (see \cite{Con00}*{Thm.~5.1.2}). If $m > 3$, then $\sX_1(m)_{\bZ_{(2)}}$ is a scheme by \cite{KM85}*{2.7.4} and \cite{Ces15a}*{4.1.4}.

\item
On the level of coarse moduli schemes, the pullback map
\[
\qq H^0(X_0(m)_{\bZ_{(2)}}, \Omega^1) \ra H^0(X_0(m)_{\bF_2}, \Omega^1)
\]
is surjective as in the proof of \ref{SC-a} (see \S\ref{str-odd-n} for basic properties of $X_0(m)_{\bZ_{(2)}}$). Therefore, the `if and only if' claim follows from Remark \ref{AM-away}, which supplies the identification 
\[
\qq H^0(X_0(m)_{\bZ_{(2)}}, \Omega^1) \cong H^0(\sX_0(m)_{\bZ_{(2)}}, \Omega^1).
\]
Granted that we address the case when there exists a suitable $q$, if $m$ is a prime, then we may assume that $m \ge 5$, so that the surjectivity of \eqref{coarse-thing} results from \cite{Maz77}*{II.4.4~(1)}.

For the rest of the proof we set $\sZ \ce \sX_0(m)_{\bF_{2}}$ and $Z \ce X_0(m)_{\bF_{2}}$ for brevity and recall that $Z$ is the coarse moduli space of $\sZ$ because $2 \nmid m$ (see \cite{Ces15a}*{6.4~(b)}). It suffices to show that
\be \lab{omega-map}
\qq \Omega^1_{Z/\bF_2} \ra \pi_* \Omega^1_{\sZ/\bF_2}
\ee
induced by pullback along the coarse moduli scheme morphism $\pi \colon \sZ \ra Z$ is an isomorphism under the assumption of the existence of $q$. The proof of this is similar to the proof of \Cref{app-main}~\ref{AM-a}, and the role of $q$ is to ensure that the ramification of $\pi$ is tame.

For every odd $m$, \eqref{omega-map} is an isomorphism over the open $V \subset Z$ on which the $j$-invariant satisfies $j \neq 0$ because $\pi|_{\pi\i(V)}$ is \'{e}tale (see \cite{Ces15a}*{proof of Thm.~6.7}) so that \eqref{omega-map} over $V$  identifies with the $\Omega^1_{V/\bF_2}$-twist of the isomorphism $\sO_{V} \isomto (\pi|_{\pi\i(V)})_*(\sO_{\pi\i(V)})$. It remains to analyze \eqref{omega-map} after base change to the completion $\wh{\sO}_{Z, z}^\sh$ of the strict Henselization of $Z$ at a variable $z \in Z(\ov{\bF}_2)$ with $j(z) = 0$. The $\bF_2$-smoothness of $Z$ and $\sZ$ gives an isomorphism
\[
\qqq \wh{\sO}_{Z, z}^\sh \simeq \ov{\bF}_2\llb t \rrb \qq \text{under which} \qq (\Omega^1_{Z/\bF_2})_{\wh{\sO}_{Z, z}^\sh} \simeq \ov{\bF}_2\llb t \rrb \cdot dt
\]
and also, using the identification $Z(\ov{\bF}_2) \cong \sZ(\ov{\bF}_2)$ to view $z$ inside $\sZ(\ov{\bF}_2)$, an isomorphism
\[
\qqq \wh{\sO}_{\sZ, z}^\sh \simeq \ov{\bF}_2\llb \tau \rrb \qq \text{under which} \qq (\Omega^1_{\sZ/\bF_2})_{\wh{\sO}_{\sZ, z}^\sh}\simeq \ov{\bF}_2\llb \tau \rrb \cdot d\tau.
\]
Moreover, with $G \ce \Aut(z)/\{\pm 1\}$ we have compatible identifications
\[
\qqq \ov{\bF}_2\llb t \rrb \cong (\ov{\bF}_2\llb \tau \rrb)^G  \qq \text{and} \qq (\pi_* \Omega^1_{\sZ/\bF_2})_{\wh{\sO}_{Z, z}^\sh} \cong (\ov{\bF}_2\llb \tau \rrb \cdot d\tau)^G
\]
with $G$ acting faithfully on $\ov{\bF}_2\llb \tau \rrb$ (see \cite{DR73}*{I.8.2.1} or \cite{Ols06}*{2.12}). If $\pi$ is tamely ramified at $z$ (i.e.,~if $2 \nmid \#G$), then $G \simeq \mu_{\#G}(\ov{\bF}_2)$ and we may choose $\tau$ in such a way that $t = \tau^{\#G}$ and any $\zeta \in \mu_{\#G}(\ov{\bF}_2)$ acts by $\tau \mapsto \zeta \cdot \tau$. 
Therefore, in the tamely ramified case the map $\ov{\bF}_2\llb t \rrb \cdot dt \ra (\ov{\bF}_2\llb \tau \rrb \cdot d\tau)^G$ that identifies with the $\wh{\sO}_{Z, z}^\sh$-pullback of \eqref{omega-map} is an isomorphism.

To complete the proof we show that $\pi$ is tamely ramified at $z$ if some prime $q$ with $q \equiv 3 \bmod 4$ divides $m$. Let $E \ra \Spec \ov{\bF}_2$ be an elliptic curve that underlies $z \in \sZ(\ov{\bF}_2)$. The action of $\Aut(z)$ on $E[q](\ov{\bF}_2)$ is faithful (because $q \ge 3$) and preserves the Weil pairing and a cyclic subgroup $C$ of order $q$. Thus, since a $2$-Sylow subgroup of $\Aut(z)$ acts semisimply, its action on $C$ embeds it into $\Aut(C) \simeq (\bZ/q\bZ)^\times$. To conclude that the inclusion of $\{\pm 1\}$ into this $2$-Sylow subgroup is an equality, as desired, it remains to note that $\#((\bZ/q\bZ)^\times[2^\infty]) = 2$ because $q \equiv 3 \bmod 4$.
\qedhere
\eenum
\epf

\bremst 
\remit
The equivalent conditions of \Cref{surj-crit}~\ref{SC-b} are also equivalent to the inequality 
\[
\qq \dim_{\bF_2} H^0(\sX_0(m)_{\bF_2}, \Omega^1) \le g, \qq \text{where $g$ is the genus of $X_0(m)_\bC$.}
\]
To see this it suffices to note that $\dim_{\bF_2} H^0(X_0(m)_{\bF_2}, \Omega^1) = g$ and that the generic isomorphy of the map \eqref{omega-map} ensures the injectivity of \eqref{coarse-thing}. 

\remit \lab{always}
We are unaware of examples of odd $m$ for which the equivalent conditions of \Cref{surj-crit}~\ref{SC-b} fail to hold. The proof of \Cref{manin-main}~\ref{MM-1}--\ref{MM-15} given in \S\S\ref{stevens-pf-1}--\ref{manin-pf-1} shows that if these conditions hold for $m = \f{n}{2}$, then the Manin constant of any new elliptic optimal quotient of $J_0(n)$ is odd.
\eremst

%% file: proofs.tex

\section{Proofs of the main results}

With the results of \S\S\ref{reduction}--\ref{oldforms} at our possession, we are ready to present the proofs of \Cref{stevens-main,manin-main}. Most of the $p = 2$ cases of these theorems are proved in \S\S\ref{stevens-pf-1}--\ref{manin-pf-1}, whereas the remaining cases are postponed until \S\S\ref{stevens-pf-2}--\ref{manin-pf-2} because they rely on a direct relationship between the conjectures of Manin and Stevens, a relationship captured by \Cref{jj-ee} and encapsulated by the formula~\eqref{MS-reln}.

\bpp[Proof of \Cref{stevens-main} in the case $p = 2$] \lab{stevens-pf-1}
For a new elliptic optimal quotient
\[
\pi\colon J_1(n) \surjects E \qq \text{with} \qq \ord_2(n) \le 1,
\]
we seek to show that $\ord_2(c_\pi) = 0$. \Cref{odd-n} settles the case of an odd $n$, so we assume that $\ord_2(n) = 1$. In this case, since $c_\pi \in \bZ$ (see \S\ref{setup}),  \Cref{main-reduction} reduces us to showing that
\be \lab{seek-0}
\tst \f{H^0(U^\mu, \Omega^1)}{H^0(X_1(n)_{\bZ_{(2)}}, \Omega) + H^0(U^\mu, \Omega^1) \cap (\bQ \cdot f)^\perp} = 0,
\ee
where $f$ is the normalized new eigenform that corresponds to $\pi$ and $U^\mu \subset X_1(n)_{\bZ_{(2)}}$ is the complement of an irreducible component of $X_1(n)_{\bF_{2}}$. Since the pullback map
\[
\tst H^0(\sX_1(\f{n}{2})_{\bZ_{(2)}}, \Omega^1) \ra H^0(\sX_1(\f{n}{2})_{\bF_2}, \Omega^1)
\]
is surjective by \Cref{surj-crit}~\ref{SC-a}, \Cref{lift-all} shows that $H^0(U^\mu, \Omega^1)/H^0(X_1(n)_{\bZ_{(2)}}, \Omega)$ consists of images of oldforms in $H^0(U^\mu, \Omega^1)$. To obtain \eqref{seek-0}, it therefore remains to note that every such oldform lies in $H^0(U^\mu, \Omega^1) \cap (\bQ \cdot f)^\perp$.
\QED
\epp

\bpp[Proof of \Cref{manin-main}~\ref{MM-1}--\ref{MM-15}] \lab{manin-pf-1}
For a new elliptic optimal quotient
\[
\pi\colon J_0(n) \surjects E \qq \text{with} \qq \ord_2(n) \le 1,
\]
we seek to show that $\ord_2(c_\pi) = 0$ whenever $n$ has a prime factor $q$ with $q \equiv 3 \bmod 4$ and whenever $n = 2p$ for some prime $p$. The argument is the same as that of \S\ref{stevens-pf-1}, except that we use \Cref{surj-crit}~\ref{SC-b} in place of \Cref{surj-crit}~\ref{SC-a}.
\QED
\epp

The proof of \Cref{stevens-main} for an odd $p$ is given in \S\ref{stevens-pf-2} and proceeds by reduction to $J_0(n)$, more precisely, to \Cref{odd-p}. Even though we have already proved this proposition in Remark \ref{odd-reproof}, we also include its more standard proof based on exactness properties of semiabelian N\'{e}ron models.

\bpropt \lab{odd-p}
For a new elliptic optimal quotient $\pi\colon J_0(n) \surjects E$ and an odd prime $p$ with $\ord_p(n) \le 1$, the Manin constant $c_\pi \in \bZ$ satisfies $\ord_p(c_\pi) = 0$.
\epropt

\bpf
Since $J_0(n)$ has semiabelian reduction at $p$ and $p - 1 > 1$, \cite{BLR90}*{7.5/4 and its proof} ensure that $\pi$ induces a smooth map $\pi_{\bZ_{(p)}}\colon \cJ_{\bZ_{(p)}} \ra \cE_{\bZ_{(p)}}$ on the N\'{e}ron models over $\bZ_{(p)}$. Due to the resulting surjectivity of $\Lie(\pi_{\bZ_{(p)}})\colon \Lie \cJ_{\bZ_{(p)}} \ra \Lie \cE_{\bZ_{(p)}}$,  the dual map $\iota\colon H^0(\cE_{\bZ_{(p)}}, \Omega^1) \hra H^0(\cJ_{\bZ_{(p)}}, \Omega^1)$ has a torsion free cokernel. Since $\im(\iota) = \bZ_{(p)} \cdot c_\pi f$, where $f$ is the normalized newform that corresponds to $\pi$, to conclude that $\ord_p(c_\pi) = 0$ it remains to recall that $f \in H^0(\cJ_{\bZ_{(p)}}, \Omega^1)$ (see Remark \ref{also-odd-p}).
\epf

The following lemma supplies a direct relationship between \Cref{stevens-conj,manin-conj}.

\blemt \lab{jj-ee}
Let $\pi_0\colon J_0(n) \surjects E_0$ and $\pi_1\colon J_1(n) \surjects E_1$ be new elliptic optimal quotients  that correspond to the same normalized new eigenform $f$.
\benum
\item \lab{JE-a}
There is a unique isogeny $e$ that fits into the commutative diagram
\be\ba \lab{jj-diag}
\xymatrix{
J_1(n) \ar@{->>}[r]^{\pi_1} \ar@{->>}[d]^{j^\vee} & E_1 \ar@{->>}[d]^-{e} \\
J_0(n) \ar@{->>}[r]^{\pi_0} & E_0
}
\ea\ee
in which $j^\vee$ is the dual of the pullback map $j\colon J_0(n) \ra J_1(n)$. Moreover, the $\bQ$-group scheme $\Ker e$ is constant and is a quotient of the Cartier dual $\Sigma(n)^\vee$ of the Shimura subgroup 
\[
\q \Sigma(n) \ce \Ker\p{J_0(n) \xra{j} J_1(n)}.
\]

\item \lab{JE-b}
The Manin constant $c_{\pi_0}$ of $\pi_0$ and the Manin--Stevens constant $c_{\pi_1}$ of $\pi_1$ are related by
\be \lab{MS-reln}
c_{\pi_0} = c_{\pi_1} \cdot \#\Coker\p{H^0(\cE_0, \Omega^1) \xra{e^*} H^0(\cE_1, \Omega^1)},
\ee
where $\cE_0$ and $\cE_1$ are the N\'{e}ron models of $E_0$ and $E_1$ over $\bZ$. In particular, Conjecture \uref{manin-conj} for $\pi_0$ implies Conjecture \uref{stevens-conj} for $\pi_1$.

\eenum  
\elemt

\bpf \hfill
\benum
\item
The existence of the unique $e$ follows from the Hecke equivariance of $\pi_0 \circ j^\vee$. By \cite{LO91}*{Thm.~2}, the finite $\bQ$-group $\Sigma(n)$ is of multiplicative type, so $\Sigma(n)^\vee$ is constant. Thus, it suffices to argue that $\Ker e$ is a quotient of $\Sigma(n)^\vee$ or, since $\Ker \pi_0$ is connected, that the component group of $\Ker(j^\vee)$ is $\Sigma(n)^\vee$. The latter follows from the exact sequences
\[
\qq 0 \ra (\Coker j)^\vee \ra J_1(n) \xra{a} (\im j)^\vee \ra 0 \q \text{and} \q 0 \ra \Sigma(n)^\vee \ra (\im{j})^\vee \xra{b} J_0(n) \ra 0
\]
in which $\Coker j$ and $\im j$ are abelian varieties and $b \circ a = j^\vee$.

\item
The formula \eqref{MS-reln} follows from \eqref{jj-diag} once we note that the alternative description of \eqref{id-omg} reviewed in \S\ref{conv} ensures that the $j^\vee$-pullback of $f$ is $f$ (see Remark \ref{also-odd-p}). The last sentence follows from \eqref{MS-reln} because $c_{\pi_1} \in \bZ$ (see \S\ref{setup}). \qedhere
\eenum
\epf

\bpp[Proof of \Cref{stevens-main} in the case of an odd $p$] \lab{stevens-pf-2}
For an odd prime $p$ and a new elliptic optimal quotient
\[
\pi_1\colon J_1(n) \surjects E_1 \qq \text{with} \qq \ord_p(n) \le 1,
\]
we seek to show that $\ord_p(c_{\pi_1}) = 0$. For this, due to \eqref{MS-reln}, it suffices to recall that $c_{\pi_1} \in \bZ$ and to note that, by \Cref{odd-p}, the elliptic new optimal quotient $\pi_0 \colon J_0(n) \surjects E_0$ that corresponds to the same normalized new eigenform as $\pi_1$ has $\ord_p(c_{\pi_0}) = 0$.
 \QED
\epp

\bpp[Proof of \Cref{manin-main}~\ref{MM-2}] \lab{manin-pf-2}
For a new elliptic optimal quotient
\[
\pi_0\colon J_0(n) \surjects E_0 \qq \text{with} \qq \ord_2(n) \le 1,
\]
we seek to show that $\ord_2(c_{\pi_0}) = 0$ whenever $E_0(\bQ)[2] = 0$. For this, in the notation of \Cref{jj-ee} and thanks to \eqref{MS-reln}, it suffices to prove $2 \nmid \# \Ker e$ because $\ord_2(c_{\pi_1}) = 0$ by \Cref{stevens-main}. However, $\Ker e$ is constant, so it remains to note that $E_1$, being isogenous to $E_0$, satisfies $E_1(\bQ)[2] = 0$.
\QED
\epp

\bremt
By \cite{LO91}*{Cor.~2 on p.~173}, if $n = 2 \cdot q^r$ for a prime $q$ with either $q\equiv 3 \bmod 4$ or $q \equiv 5 \bmod 8$, then the order $\#\Sigma(n)$ of the Shimura subgroup is odd. Therefore, since \Cref{jj-ee}~\ref{JE-a} ensures that $\Ker e$ is a quotient of $\Sigma(n)^\vee$, for such $n$ the argument of \S\ref{manin-pf-2} shows that $\ord_2(c_{\pi_0}) = 0$ for every new elliptic optimal quotient $\pi_0\colon J_0(n) \surjects E_0$.
\eremt

%% file: app-integral.tex

\section{The ``relative dualizing sheaf'' of $\sX_H$ and of its coarse space} \lab{append}

The main goal of this appendix is to prove a certain comparison result between the relative dualizing sheaf on the modular curve $X_H$ and an analogous sheaf on $\sX_H$. This is accomplished in \Cref{app-main} after introducing the relevant sheaf on $\sX_H$ in \S\ref{rel-dual} and detailing some of its properties in \S\ref{modular-case}. As the proof of \Cref{int-thy} illustrates, the practical role of \Cref{app-main} is to facilitate passage between $X_H$ and $\sX_H$ in the study of integral structures on spaces of weight $2$ cusp forms (a link between the latter and the relative dualizing sheaf is supplied by Kodaira--Spencer, see \Cref{KS}).

\bpp[``Relative dualizing sheaves'' of Deligne--Mumford stacks] \lab{rel-dual}
For a scheme $S$ and a Cohen--Macaulay (and hence flat) morphism $X \ra S$ that has a pure relative dimension, the theory of Grothendieck duality associates a quasi-coherent, locally finitely presented, $S$-flat relative dualizing $\sO_X$-module $\Omega_{X/S}$ (see \cite{Con00}*{bottom halves of p.~157 and p.~214}), which identifies with the determinant of $\Omega^1_{X/S}$ if $X \ra S$ is in addition smooth. The formation of $\Omega_{X/S}$ is compatible with \'{e}tale localization on $X$: for every \'{e}tale $S$-morphism $f\colon X' \ra X$ one has a canonical isomorphism
\be \lab{et-omega}
\iota_f \colon f^*(\Omega_{X/S}) \isomto \Omega_{X'/S}
\ee
supplied by \cite{Con00}*{Thm.~4.3.3 and bottom half of p.~214}. Moreover, if $f' \colon X'' \ra X'$ is a further \'{e}tale $S$-morphism, then \cite{Con00}*{(4.3.7) and bottom half of p.~214} supply the following compatibility between the isomorphisms of \eqref{et-omega}:
\be \lab{iso-comp}
\iota_{f \circ f'} = \iota_{f'} \circ ((f')^*(\iota_{f})) \colon (f')^*(f^*(\Omega_{X/S})) \isomto \Omega_{X''/S}.
\ee
Therefore, if $\sX$ is a Deligne--Mumford stack over $S$ such that $\sX \ra S$ Cohen--Macaulay and has a pure relative dimension, then the compatibilities \eqref{iso-comp} ensure that the $\sO_X$-modules $\Omega_{X/S}$ for \'{e}tale morphisms $X \ra \sX$ from a scheme $X$ glue to form a quasi-coherent, locally of finite presentation, $S$-flat $\sO_{\sX}$-module 
\[
\Omega_{\sX/S}, \qq \text{the ``relative dualizing sheaf'' of} \qq \sX \ra S
\]
(see \cite{LMB00}*{12.2.1} for a discussion of analogous compatibilities). If $\sX \ra S$ is in addition smooth, then $\Omega_{\sX/S}$ identifies with the determinant of $\Omega^1_{\sX/S}$. Due to \cite{Con00}*{Thm.~4.4.4 and bottom half of p.~214}, the formation of $\Omega_{\sX/S}$ commutes with base change in $S$.
\epp

\bremt
In the case when $\sX \ra S$ is proper (and $\sX$ is not a scheme), we \emph{do not} claim any dualizing properties of the $\sO_\sX$-module $\Omega_{\sX/S}$ constructed in \S\ref{rel-dual}. Nevertheless, if a sufficiently robust Grothendieck--Serre duality formalism for $\sX_0(m)_{\bZ_{(2)}} \ra \Spec \bZ_{(2)}$ with $2\nmid m$ existed, then it would prove the surjectivity in \Cref{surj-crit}~\ref{SC-b} without assuming the primality of $m$ or the existence of $q$ (see the proof of \Cref{surj-crit}~\ref{SC-a}), which would settle the semistable case of the Manin conjecture (see Remark \ref{always}).
\eremt

\bpp[The case of modular curves] \lab{modular-case}
For us, the key case in \S\ref{rel-dual} is when $S = \Spec \bZ$ and $\sX$ is either a modular stack $\sX_H$ or its coarse moduli scheme $X_H$ for some open subgroup $H \subset \GL_2(\wh{\bZ})$ (see \S\ref{notation}), as we now assume. The resulting $\sX \ra S$ is flat, of finite presentation, purely of relative dimension $1$, and Cohen--Macaulay (the latter due to the normality of $\sX$ and \cite{EGAIV2}*{6.3.5~(i)}), so the discussion of \S\ref{rel-dual} applies. Moreover, \cite{EGAIV2}*{6.12.6 (i)} and the normality of $\sX$ ensure that after removing finitely many closed points $\sX$ becomes regular and hence also a local complete intersection over $\bZ$ (see \cite{Liu02}*{6.3.18}). In particular, each of the finitely many nonsmooth $\bZ$-fibers of $\sX$ has a dense open Gorenstein locus. The resulting coherent $\bZ$-flat $\sO_\sX$-module $\Omega_{\sX/\bZ}$ is therefore a line bundle on a $\bZ$-fiberwise dense open of $\sX$. It then follows from \cite{EGAIV2}*{6.4.1~(ii)} and from the proof of \cite{Con00}*{Lem.~5.2.1} (carried out for the compactification of an \'{e}tale scheme cover of a $\bZ$-fiber of $\sX$) that $\Omega_{\sX/\bZ}$ is Cohen--Macaulay\footnote{The $S$-fibral Cohen--Macaulayness of $\Omega_{\sX/S}$ is actually a general fact.} 
 on the entire $\sX$.
\epp

With the discussion of \S\ref{modular-case}, we are ready for the promised comparison result.

\bthmt \lab{app-main}
Fix an open subgroup $H \subset \GL_2(\wh{\bZ})$ and let $\pi\colon \sX_H \ra X_H$ be the coarse moduli space morphism. 

\benum
\item \lab{AM-a}
Pullback of K\"{a}hler differentials along $\pi_\bQ$ induces an $\sO_{(X_H)_\bQ}$-module isomorphism
\be \lab{omega-comp}
\q \Omega^1_{(X_H)_\bQ/\bQ} \isomto (\pi_\bQ)_*( \Omega^1_{(\sX_H)_\bQ/\bQ}).
\ee

\item \lab{AM-b}
For every open subscheme $U \subset X_H$ with $\sU \ce \pi\i(U)$ such that $\pi|_{\sU}\colon \sU \ra U$ is \'{e}tale over a $\bZ$-fiberwise dense open of $U$, the isomorphism $H^0(U_\bQ, \Omega^1) \cong H^0(\sU_\bQ, \Omega^1)$ of \ref{AM-a} identifies
\[
\qq H^0(U, \Omega) \subset H^0(U_\bQ, \Omega^1) \qq \text{with} \qq H^0(\sU, \Omega) \subset H^0(\sU_\bQ, \Omega^1).
\]
\eenum
\ethmt

\bpf \hfill
\benum
\item
Let $V \subset (X_H)_\bQ$ be a dense open over which $\pi_\bQ$ is \'{e}tale, and set $\sV \ce (\pi_\bQ)\i(V)$. The  restriction of \eqref{omega-comp} to $V$ identifies with the $\Omega^1_{V/\bQ}$-twist of the isomorphism $\sO_V \isomto (\pi_\bQ|_\sV)_*(\sO_\sV)$, so is an isomorphism. It remains to prove that the base change of \eqref{omega-comp} to the completion $\wh{\sO}_{(X_H)_\bQ, x}^\sh$ of the strict Henselization of $(X_H)_\bQ$ at a variable $x \in X_H(\ov{\bQ})$ is an isomorphism.

We have an isomorphism 
\[
\qq \wh{\sO}_{(X_H)_\bQ, x}^\sh \simeq \ov{\bQ}\llb t\rrb \qq \text{under which} \qq (\Omega^1_{(X_H)_\bQ/\bQ})_{\wh{\sO}_{(X_H)_\bQ, x}^\sh} \simeq \ov{\bQ}\llb t\rrb \cdot dt,
\]
and also (using the identification $X_H(\ov{\bQ}) \cong \sX_H(\ov{\bQ})$ to view $x$ inside $\sX_H(\ov{\bQ})$)
\[
\qq \wh{\sO}_{(\sX_H)_\bQ, x}^\sh \simeq \ov{\bQ}\llb \tau\rrb \qq \text{under which} \qq (\Omega^1_{(\sX_H)_\bQ/\bQ})_{\wh{\sO}_{(\sX_H)_\bQ, x}^\sh} \simeq \ov{\bQ}\llb \tau\rrb \cdot d\tau.
\]
Taking into account the action of the automorphism group of $x \in \sX_H(\ov{\bQ})$, we have, compatibly, 
\[
\qqq\wh{\sO}_{(X_H)_\bQ, x}^\sh \cong (\wh{\sO}_{(\sX_H)_\bQ, x}^\sh)^G \q \text{and} \q ((\pi_\bQ)_*( \Omega^1_{(\sX_H)_\bQ/\bQ}))_{\wh{\sO}_{(X_H)_\bQ, x}^\sh} \cong ((\Omega^1_{(\sX_H)_\bQ/\bQ})_{\wh{\sO}_{(\sX_H)_\bQ, x}^\sh})^G
\]
for a certain group $G$ acting faithfully on $\wh{\sO}_{(\sX_H)_\bQ, x}^\sh$ (see \cite{DR73}*{I.8.2.1} or \cite{Ols06}*{2.12}). Moreover, the ramification of $\pi_\bQ$ is tame, so we may assume that $G \simeq \mu_{\#G}(\ov{\bQ})$ and that $\zeta \in \mu_{\#G}(\ov{\bQ})$ acts by $\tau \mapsto \zeta \cdot \tau$ with $t = \tau^{\#G}$. It then follows that $\ov{\bQ}\llb t\rrb \cdot dt \isomto (\ov{\bQ}\llb \tau\rrb \cdot d\tau)^G$, i.e., that the base change of \eqref{omega-comp} to $\wh{\sO}_{(X_H)_\bQ, x}^\sh$ is an isomorphism. 

\item
Let $U' \subset U$ be a $\bZ$-fiberwise dense open over which $\pi$ is \'{e}tale and let $\sU' \subset \sU$ be its preimage. By \S\ref{modular-case}, the $\sO_{X_H}$-module $\Omega_{X_H/\bZ}$ is (S$_2$), and likewise for $\Omega_{\sX_H/\bZ}$, so, due to \cite{EGAIV2}*{5.10.5},
\[ \ba
\qq H^0(U, \Omega) = H^0(U', \Omega) \cap H^0(U_\bQ, \Omega^1) \qq &\text{inside} \q H^0(U'_\bQ, \Omega^1), \qq \text{and} \\
\qq H^0(\sU, \Omega) = H^0(\sU', \Omega) \cap H^0(\sU_\bQ, \Omega^1) \qq &\text{inside} \q H^0(\sU'_\bQ, \Omega^1).
\ea \] 
Therefore, \ref{AM-a} reduces us to the case when $U = U'$. Moreover, the (S$_2$) property ensures that neither $H^0(U, \Omega)$ nor $H^0(\sU, \Omega)$ changes if we remove finitely many closed points from $U$, so, thanks to \S\ref{modular-case}, we assume further that $U$ and $\sU$ are regular and that $\Omega_{U/\bZ}$ and $\Omega_{\sU/\bZ}$ are line bundles. Then, due to the \'{e}taleness of $\pi|_\sU$, we have $(\pi|_{\sU})^* (\Omega_{U/\bZ}) \cong \Omega_{\sU/\bZ}$ (see \eqref{et-omega}), so that, since $\Omega_{U/\bZ}$ is a line bundle, the resulting pullback map
\be \lab{omg-U-pullb}
\qq \Omega_{U/\bZ} \ra (\pi|_\sU)_*(\Omega_{\sU/\bZ})
\ee
identifies with the $\Omega_{U/\bZ}$-twist of the isomorphism $\sO_U \isomto (\pi|_\sU)_*(\sO_\sU)$ and hence is an isomorphism. The sought claim then follows by taking global sections in \eqref{omg-U-pullb}. \qedhere
\eenum
\epf

\bremst
\remit \lab{AM-away}
If $H$ contains $\Ker(\GL_2(\wh{\bZ}) \surjects \GL_2(\bZ/n\bZ))$, then the $\bZ$-fibral generic \'{e}taleness assumption of \Cref{app-main}~\ref{AM-b} holds for every $U$ on which $n$ is invertible (see \cite{Ces15a}*{last paragraph of the proof of Prop.~6.4~(b)}). In particular, since $\sX_H$ and $X_H$ are $\bZ[\f{1}{n}]$-smooth (see \cite{DR73}*{IV.6.7 and VI.6.7} or also \cite{Ces15a}*{6.4~(a)}), \Cref{app-main} proves that pullback induces an isomorphism
\[
\qq \Omega^1_{(X_H)_{\bZ[\f{1}{n}]}/\bZ[\f{1}{n}]} \isomto (\pi_{\bZ[\f{1}{n}]})_*(\Omega^1_{(\sX_H)_{\bZ[\f{1}{n}]}/\bZ[\f{1}{n}]}).
\]

\remit \lab{AM-X0n}
An important case in which the $\bZ$-fibral generic \'{e}taleness assumption of \Cref{app-main}~\ref{AM-b} holds for every $U$ is when $H = \Gamma_0(n)$ (see \cite{Ces15a}*{proof of Thm.~6.7}). In this case, \Cref{app-main} provides an $\sO_{X_0(n)}$-module isomorphism 
\[
\qq \Omega_{X_0(n)/\bZ} \isomto \pi_*(\Omega_{\sX_0(n)/\bZ})
\]
which on the $\bQ$-fiber is induced by pullback of K\"{a}hler differentials.
\eremst